

\input amstex
\documentstyle{amsppt}
\magnification = \magstep1
\hsize = 6.25 truein
\vsize = 22 truecm
\baselineskip .22in
\define\Lap{\varDelta}
\define\del{\partial}
\define\a{\alpha}
\predefine\b{\barunder}
\redefine\b{\beta}
\define\ga{\gamma}
\predefine\d{\dotunder}
\redefine\d{\delta}
\define\th{\theta}

\define\e{\epsilon}
\define\s{\sigma}
\predefine\o{\orsted}
\redefine\o{\omega}

\define\cin{\Cal C^{\infty}}
\define\RR{\Bbb R}
\define\SS{\Bbb S}
\define\z{\zeta}
\define\Sn{{\Bbb S}^n}

\define\SL{\Sn \backslash \Lambda}
\define\ML{\Cal M_{\Lambda}}
\define\CML{\overline{\Cal M_{\Lambda}}}
\define\ind{\text{ind\,}}
\define\relind{\text{rel-ind\,}}
\define\MML{\Bbb M_{\Lambda}}
\define\BL{\Cal B}
\define\tr{\text{tr\,}}


\NoRunningHeads
\topmatter
\title Moduli Spaces of Singular Yamabe Metrics \endtitle
\author Rafe Mazzeo${}^{(\dag)}$, Daniel Pollack${}^{(\ddag)}$ and
Karen Uhlenbeck${}^{(\star)}$ \endauthor
\abstract
Complete, conformally flat metrics of constant positive scalar curvature on
the complement of $k$ points in the $n$-sphere, $k \ge 2$, $n \ge 3$,
were constructed by R\. Schoen [S2]. We consider the problem
of determining the moduli space of all such metrics.
All such metrics are asymptotically
periodic, and we develop the  linear analysis necessary
to understand the nonlinear problem. This includes a Fredholm
theory and asymptotic regularity theory for the Laplacian
on asymptotically periodic manifolds, which is of
independent interest. The main result is that
the moduli space is a locally real analytic variety of dimension
$k$. For a generic set of nearby conformal classes the moduli  space is
shown to be a $k-$dimensional real analytic manifold.
The structure as a real
analytic variety is obtained by writing the space as an intersection
of a Fredholm pair of infinite dimensional real analytic manifolds.
\endabstract
\affil Mathematical Sciences Research Institute ($\dag$ and $\ddag$),
Stanford University ($\dag$) and University of Texas at Austin
($\ddag$ and $\star$) \endaffil
\thanks Research supported in part by ($\dag$) NSF Young Investigator
Award, the Sloan Foundation, NSF grant \# DMS9001702 and  ($\ddag$)
NSF grant \# DMS9022140.
\endthanks
\endtopmatter
\document

\specialhead I. Introduction \endspecialhead

Much has been clarified in the past ten years about the behavior
of solutions of the semilinear elliptic equation relating the scalar
curvature functions of two conformally related metrics. The starting
point for these recent developments was R.~Schoen's
resolution of the Yamabe problem on compact
manifolds [S1], capping the work of a number of mathematicians
over many years.  Soon thereafter Schoen [S2] and Schoen-Yau [SY]
made further strides in understanding weak solutions of this
equation, particularly on the sphere, and its relationship with
conformal geometry. Of particular interest here is the former,
[S2]; in that paper, Schoen constructs metrics with
constant positive scalar curvature on $\SS^n$, conformal to the standard
round metric, and with prescribed isolated singularities (he also
constructs solutions with certain, more general, singular sets).
This `singular Yamabe problem' is to find a metric $g = u^{4/(n-2)}g_0$
on a domain $\SS^n \backslash \Lambda$ which is complete and has constant
scalar curvature $R(g)$. This is equivalent to finding a positive
function $u$ satisfying
$$
\gathered \Lap_{g_0} u  - \frac{n-2}{4(n-1)}R(g_0)u +
\frac{n-2}{4(n-1)}R(g)u^{\frac{n+2}{n-2}} = 0 \ \text{on}\ \SS^n \backslash
\Lambda \\ \text{$g$ complete on $\SS^n \backslash \Lambda$},
\qquad R(g) = \ \text{constant,} \endgathered \tag 1.1
$$
where $R(g_0) = n(n-1)$ is the scalar curvature of the round metric $g_0$.
The completeness of $g$ requires that $u$ tend to infinity, in
an averaged sense, on approach to $\Lambda$.

The earliest work on this singular Yamabe problem seems to have
been that of C.~Loewner and L.~Nirenberg [LN], where metrics with
constant negative scalar curvature are constructed. Later
work on this `negative' case was done by Aviles-McOwen [AMc]
cf\. also [Mc] and Finn-McOwen [FMc], where the background
manifold and metric are allowed to be arbitrary. For a solution
with $R(g) < 0$ to exist, it is necessary and sufficient that $\dim(
\Lambda) > (n-2)/2$.  A partial converse is that if a solution
with $R(g) \ge 0$ exists (at least when $M = \SS^n$), it is necessary
that $k \le (n-2)/2$, [SY]. Schoen gave the first general construction
of solutions with $R(g) > 0$ [S2]. These solutions have singular set
$\Lambda$ which is either discrete or nonrectifiable.  Many new solutions
on the sphere with singular set $\Lambda$,
a smooth perturbation of an equatorial $k$-sphere,
$1\le k < (n-2)/2$, are constructed in [MS].
F.~Pacard [Pa] has recently constructed positive complete
solutions on the sphere with singular set an arbitrary smooth
submanifold of dimension $(n-2)/2$.  Regularity
of solutions near the singular set $\Lambda$, no matter the sign
of the curvature or background manifold or metric, is
examined in [M1], cf\. also [ACF] for a special case of
relevance to general relativity.  A more detailed account of
part of this history is given in [MS].

In this paper we return to this problem in the setting studied by
Schoen [S2] on $\SS^n$, where $\Lambda$ is a finite point-set:
$$ \Lambda = \{p_1, \cdots, p_k\}.  \tag 1.2 $$
Hereafter, $\Lambda$ will always be taken to be this set, unless
indicated specifically, and the scalar curvature $R(g)$ attained
by the conformal metric will always be
$n(n-1)$. In this case, no solutions of (1.1) exist when $k=1$.
A proof of this, following from a general symmetry theorem, is indicated
below.  Thus we always assume that $k \ge 2$.
Since this problem is conformally invariant, the
set $\Lambda$ may be replaced by $F(\Lambda)$ for any conformal
transformation $F \in O(n+1,1)$. A simple topological argument
(given in [S2]) shows that using such a transformation we can
always arrange that $\Lambda$ is `balanced,' i\.e\. that
the points $p_j$ sum to zero as vectors in $\RR^{n+1}$.
Henceforth this will be tacitly assumed as well.
Notice that if $\Lambda$ is balanced, $k > 2$, and if $F$ is
a conformal transformation which preserves $\Lambda$, then (since some
power of $F$ fixes three points on the sphere, and also the origin
in $\RR^{n+1}$) $F$ must be an orthogonal transformation.
If $k=2$ then $F$ could also be a dilation. Also, any balanced configuration
$\Lambda$ is contained in a minimal equatorial subsphere $\SS^k \subset
\SS^n$, and an Alexandrov reflection argument similar to the one
in \S8 of [CGS] shows that if
$u$ is an arbitrary solution to (1.1), and $F$ is a rotation
preserving $\Lambda$ pointwise, then $u$ must also be invariant under $F$.
Observe that this implies that no complete solutions
exist on $\Sn \backslash \{p\}$. If $u$ were such a solution,
then because any other point can be moved to be antipodal to $p$, this
reflection argument would show that $u$ is rotationally
symmetric with respect to any other point on the sphere.
Hence $u$ would have to be constant, which contradicts completeness.
Hence (except when $k=2$) we do not get `trivial' families of solutions
of (1.1) obtained by pulling back a fixed solution by a family
of conformal transformations. Note, however, that if $F$ preserves
$\Lambda$, but permutes the points, there is no reason to expect
that most solutions of (1.1) will be fixed by $F$; indeed, it is
not even clear a priori that $F$-invariant solutions exist in this case.

Our aim in this work is to consider the moduli space $\ML$, which is
by definition the set of all smooth positive functions $u$ to the problem
$$
\gathered \Lap_{g_0} u  - \frac{n(n-2)}{4}u +
\frac{n(n-2)}{4}u^{\frac{n+2}{n-2}} = 0 \ \text{on}\ \SS^n \backslash
\Lambda \\
\text{$g=u^{\frac{4}{n-2}}g_{0}$ complete on $\SS^n \backslash \Lambda$},
 \endgathered \tag 1.3
$$
where $\Lambda = \{p_1, \cdots, p_k\}$ is any fixed, balanced set of
$k$ points in $\SS^{n}$.  We remark that the geometric condition that the
metric $g$ be complete is equivalent to the analytic one of requiring
that the set $\Lambda$ consist of nonremovable singularities of $u$.
This space might be called the `PDE moduli
space' $\ML^{\text{PDE}}$, to distinguish it
from the `geometric moduli space' $\ML^{\text{geom}}$ which consists of
all geometrically distinct solutions of this problem. In most instances
these spaces coincide, although the second could conceivably be somewhat
smaller than the first if some solutions admit nontrivial isometries
(e\.g\. as in the case $k=2$ discussed in \S2 below). $\ML$ will always
denote the former of these spaces here. By Schoen's work [S2],
$\ML$ is nonempty whenever $k>1$. In fact, his construction
yields families of solutions. As this was not the aim
of his work, this is not made explicit there, nor are the
free parameters in his construction counted.

We shall examine a number of questions, both local and global, concerning
the nature of this moduli space. The simplest of these is whether
$\ML$ is a manifold, or otherwise tractable set, and if so, what is its
dimension. Our main result is
\proclaim{Theorem 1.4} $\ML$ is locally a real analytic
variety of formal dimension $k$.
\endproclaim
The formal dimension is the dimension predicted by an index
theorem. As is usual in moduli space theories, obstructions
may well exist to prevent $\ML$ from attaining this dimension.
We clarify this and give a more careful statement
in Theorem 5.4,  Corollary 5.5 and Theorem 6.13 below.
We also describe how natural parameters
on $\ML$ may arise. On the linear level, these are given by
certain scattering theoretic information for the metrics $g \in \ML$.
Another geometric description is given by the Poho\v{z}aev invariants
which are defined in \S3. We also obtain information on a geometrically
natural compactification of $\ML$ which is obtained by taking the union with
lower-dimensional
moduli spaces corresponding to singular sets $\Lambda' \subset \Lambda$.
We are, as yet, unable to provide a satisfactory description of the interior
singularities of $\ML$, or to determine whether this compactification,
$\CML$, is itself a compact real (semi-)analytic variety.
This latter property is, by all indications, true.
We hope to return to this later.

It is important to mention that there are
many similarities between the theory of constant scalar
curvature metrics on $\SL$ where $\Lambda$ is finite and that of embedded,
complete, constant mean curvature (CMC) surfaces with $k$ ends in $\RR^3$.
The first examples of such CMC surfaces were given by N\.~Kapouleas [Kap].
One-parameter families of solutions with symmetry were constructed by
K\.~Grosse-Brauckmann [B].
Further general results on the structure of these surfaces and related
problems appear in
[KKS], [KKMS] and [KK]. In the last of these, N\.~Korevaar and
R\.~Kusner conjecture
that there is a good moduli space theory for these surfaces. Our methods
may be adapted immediately to this setting. If $M \subset \RR^3$ is an
embedded, complete CMC surface with $k\geq 3$ ends which satisfies
a hypothesis analogous to (5.2) below,
then our results imply that the space of all nearby surfaces of this type, up
to rigid motion in $\RR^3$, is a $(3k-6)$-dimensional real analytic manifold.
In order to understand the structure of the moduli space near surfaces
where this hypothesis is not satisfied an argument different from the one
employed here is needed since, for example, we do not have an analogue of
the constructions in \S6.
Recently such an arguement was derived based on the linear analysis
developed in this paper.  Thus we
can also assert that the moduli space of such constant mean curvature
surfaces is locally a real analytic variety, as was claimed in [KK].
This arguemnt also applies to give a new, direct proof of Theorem 1.4.
These results will appear in a forthcoming paper with R. Kusner.

The outline of this paper is as follows. In \S2 we analyze in detail
the special case when $k = 2$. The solutions here will be called
Delaunay solutions, in analogy with a similar family of
complete constant mean curvature surfaces in $\RR^3$
discovered in 1841 by C. Delaunay [D].
Only in this case may the moduli space be determined completely,
since, using the symmetry argument discussed above,
the equation now reduces to an ODE. We also analyze
the spectral theory of the linearized scalar curvature operator
completely in this simple case. \S3 collects a number of
disparate results about solutions of (1.3) with isolated singularities
which are used throughout the rest of the paper. This contains an
explanation of the
results of Caffarelli-Gidas-Spruck and Aviles-Korevaar-Schoen which
state that the Delaunay solutions are good models for
arbitrary solutions of (1.3) with isolated singularities.
We also discuss here the Poho\v zaev invariants and
compactness results for solutions $g\in \ML$.
The linearization $L$ of the scalar curvature operator
is studied in \S4. We prove Fredholm results for this operator,
using and expanding earlier work of C.~Taubes [T]. Other results here
include more detailed information on asymptotics of solutions of
$Lw=0$, as well as the computation, using a relative index
theorem, of the dimension of the `bounded
nullspace' of $L$. \S5 uses these results
to establish the structure of $\ML$ near its smooth points.
In \S6 we show that $\ML$ is a real analytic set by writing it
as a slice of an `urmoduli space' $\MML$ with the conformal
class determined by any $g \in \ML$.
Our study of the urmoduli space $\MML$ draws on work of
A.~Fischer and J.~Marsden, [FM1] and [FM2].
We also prove a generic slice theorem
here which shows that slices of $\MML$ by generic nearby conformal classes
are smooth, even if $\ML$ is not. In \S7, our concluding remarks,
we discuss three aspects about the nature of $\ML$ concerning which
we have not yet obtained satisfactory results. The first of these is the
nonexistence of $L^2$ eigenvalues for the linearization, especially
for the solutions constructed in [S2]. Secondly, we give a
description of local coordinates on $\ML$ near smooth points.
This is done both on the linear level
and geometrically in terms of the Poho\v{z}aev invariants introduced in \S3.
At this point we also provide a brief discussion of the recent
construction [MPU] of dipole solutions for the problem.
Finally, we address certain natural questions concerning the boundary of the
geometric compactification $\CML$.

The authors wish to thank Rick Schoen for his continued interest
and substantial advice throughout the course of this work.
The first two authors also had a number of enlightening conversations
with Nick Korevaar, Rob Kusner and Tom Mrowka.

\specialhead II. Delaunay Solutions \endspecialhead

In this section we discuss the Delaunay family of solutions. These
constitute, up to conformal equivalence, the totality of solutions when
$\Lambda$ has only two elements.

\head ODE analysis \endhead When $k = 2$, the condition that
$\Lambda$ be balanced means that $p_2 = -p_1$.  It is not difficult
to show using the Alexandrov reflection argument [CGS]
that any positive solution
of the PDE (1.3) is only a function of the geodesic distance from either
$p_1$ or $p_2$. This equation reduces to an ODE, which takes the simplest
form when written relative to the background metric $dt^2 + d\theta^2$ on
the cylinder $\RR \times \SS^{n-1}$ with coordinates $(t,\theta)$, which is
conformally equivalent to $\SS^n \backslash \Lambda$. Thus,
since the cylinder has scalar curvature $(n-1)(n-2)$, and $g =
u^{4/(n-2)}(dt^2+ d\theta^2)$ has scalar curvature $n(n-1)$,
$u = u(t)$ satisfies
$$ \frac {d^2}{dt^2} u - \frac {(n-2)^2}4 u + \frac {n(n-2)}4 u^{
\frac {n+2}{n-2}} = 0. \tag 2.1 $$
This is easily transformed into a first order Hamiltonian
system: setting $v = u'$ (primes denoting differentiation by $t$)
$$ \aligned
&u' = v\\ & v' = \frac{(n-2)^2}4 u - \frac{n(n-2)}4 u^{\frac {n+2}{n-2}}
\endaligned \tag 2.2 $$
The corresponding Hamiltonian energy function is
$$
H(u,v) = \frac{v^2}{2} + \frac{(n-2)^2}{8} u^{\frac {2n}{n-2}} -
\frac {(n-2)^2}{8}u^2. \tag 2.3
$$

The orbits of (2.2) remain within level sets of $H$, and since these
level sets are one dimensional, this determines these orbits (but not
their parametrizations) explicitly. The equilibrium points for this
flow are at $(0,0)$ and $(\bar u,0)$, where
$$ \bar u = \left(\frac{n-2}n\right)^{\frac {n-2}4} \tag 2.4  $$
There is a special homoclinic orbit $(u_0(t),v_0(t))$ corresponding to the
level set $H = 0$; it limits on the origin as $t$ tends to either
$\pm \infty$, and encloses a bounded set $\Omega$ in the right half-plane
which is symmetric across the $u$-axis, given by
$\{H \le 0\}$. Somewhat fortuitously
we may calculate explicitly that
$$ u_0(t) = (\cosh t)^{\frac {2-n}2}. $$
Of course, $\{H=0\}$ decomposes into two orbits: this one and the
stationary orbit $(0,0)$. It is simple to check that orbits not enclosed
by this level set, i\.e\.
those on which $H > 0$, must pass across the $v$-axis and into
the region where $u < 0$. Thus, since we are only interested in solutions
of (2.1) which remain positive and exist for all $t$, it suffices to
consider only those orbits in $\Omega$. Notice that the second
equilibrium point $(\bar u, 0)$ is in this region, and that all other
orbits are closed curves.  These correspond to periodic orbits
$(u_{\e}(t),v_{\e}(t))$, with period $T(\e)$, $0 < \e < 1$.
The parameter $\e$ may be taken as the smaller of the two
$u$ values where the orbit intersects the $u$-axis, so that
$ 0 < \e \le \bar u$ (note that, strictly speaking, the orbit
with $\e=0$ corresponds to the equilibrium point at the origin,
but by convention we set it equal to the one previously
defined).  This ODE analysis is also described in [S3].

The corresponding metrics on $\RR \times \SS^{n-1}$ (or $\SS^n \backslash
\{p_1, p_2 \}$) have a discrete group of isometries, given on the
cylinder by the translations $t \mapsto t + T(\e)$. They interpolate
between the cylindrical metric $\bar u^{4/(n-2)}(dt^2 + d\theta^2)$
(which is rescaled by the power of $\bar u$ so that its
scalar curvature is $n(n-1)$) and the solution corresponding
to the conformal factor $u_0(t)$. This later solution is nothing
other than the standard round metric on $\SS^n \backslash \{p_1,p_2\}$,
which is therefore incomplete and not, strictly speaking, in the
moduli space $\Cal M_{\Lambda}$.  In all that follows we shall
adopt the notation
$$ g_\e = u_{\e}^\frac 4{n-2}\left( dt^2 + d\theta^2\right) \tag 2.5
$$
when referring to these Delaunay metrics.

Of greatest concern is the Laplacian for the metrics $g_\e$;
in terms of the coordinates $(t,\theta)$ on the cylinder
we may write this operator as
$$ \aligned \varDelta_\e  & =   u^{-\frac {2n}{n-2}}\del_t
\left( u^2 \del_t \right) + u^{-\frac {4}{n-2}}
\varDelta_{\theta} \\
& = u^{-\frac {4}{n-2}}\del_t^2
+ 2(\del_t u)u^{-\frac {(n+2)}{n-2}}\del_t
+ u^{-\frac {4}{n-2}}\varDelta_\theta.
\endaligned \tag 2.6
$$
We use $\del_t$ for the partial derivative
with respect to $t$, etc\. Also, $\varDelta_\theta$ is the
Laplacian on the sphere $\SS^{n-1}$
of curvature $1$. We use the convention that $-\varDelta$ is
a positive operator. Also, we write $u=u_{\e}$ throughout this section.

It is convenient to replace the variable $t$ by a new variable $r(t)$,
depending on $\e$, which represents geodesic distance with respect
to $g_{\e}$ along the $\nabla t$ integral curves, which are already
geodesics (with varying parametrization) for each of the metrics $g_{\e}$.
$r(t)$ is defined by setting $dr/dt =
u^{\frac 2{n-2}}$ and $r=0$ when $t=0$. Then
$$ \del_t = u^{\frac 2{n-2}} \del_r $$
and consequently, using the first equality in (2.6),
$$\aligned
\varDelta_\e & = u^{\frac {2-2n}{n-2}}\del_r
\left( u^{\frac {2n-2}{n-2}} \del_r \right) + u^{-\frac {4}{n-2}}
\varDelta_\theta \\
& = \del_r^2 + \frac {2n-2}{n-2} \frac{\del_r u}{u}\del_r
+ u^{-\frac {4}{n-2}}\varDelta_\theta. \endaligned \tag 2.7 $$

The function $u$ is still periodic with respect to $r$, but the
period $R(\e)$ now is better behaved than the period $T(\e)$
above. In fact, whereas
$$ \lim_{\e \rightarrow 0} T(\e) = \infty,\qquad \text{and\ }
\lim_{\e \rightarrow \bar u} T(\e) = \frac {2\pi}{\sqrt {n-2}} $$
(the latter equality is proved by linearizing (2.2) at $(\bar u,0)$),
we have instead
$$ \lim_{\e \rightarrow 0} R(\e) = \pi, \qquad
\text{and } \lim_{\e \rightarrow \bar u} R(\e) = \frac {2\pi}{\sqrt n},
\tag 2.8 $$
so that the period stays within a compact interval in the positive
real axis as $\e$ varies between its limits.

It is somewhat amusing that when $n=4$ we can solve (2.1)
explicitly, once the transformation from $t$ to $r$ is
effected; this was pointed out to us by R. Schoen.  In fact,
now $\del_t = u\del_r$ and (2.1) becomes
$$
u\del_r(u\del_r u) - u + 2u^3 = 0. \tag 2.9
$$
Setting $v = u\del_r u$ as the new dependent variable and
$u$ as the independent variable we get that $ v\del_u v - u + 2u^3 = 0$.
Integrating to solve for $v$, and then integrating again leads
to a general expression for $u$. Adjusting the constants we finally
get that
$$
u_{\e}(r) = \sqrt{ \frac 12 + (\frac 12 - \e^2)\cos(2r)}. \tag 2.10
$$
The parametrization has been arranged so that $u(\e)$ attains
its minimum when $r = \pi/2 + \ell \pi, \ell \in \Bbb Z$. Notice
also that this shows that the period $R(\e) \equiv \pi$ when $n=4$.

It is perhaps also instructional, although not necessary
within the context of this paper, to consider the space of solutions
to the problem on the real projective space with one point
removed,  $\RR P^n \backslash \{p\}$. We can,
identify any solution on this space with a solution
on $\Sn \backslash \{p,-p\}$ invariant under the
reflection which exchanges these two antipodal points.
Transforming the twice-punctured sphere to the cylinder, we are looking for
solutions $u$ such that $u(t,\theta) = u(-t,-\theta)$.
Any such solution is, of course, independent of $\theta$ and in the Delaunay
family. For each value of the Delaunay parameter $\e$ there
are two possibilities, one taking its minimum value at $t=0$
and the other taking its maximum value there. These are
connected via the cylindrical solution by letting $\e$
increase to $\bar u$. The moduli space is thus a copy
of $\RR$, with the cylindrical solution at the `origin.'
Alternatively, this moduli space is just the part of
the $u$-axis in the set $H <0$ in the $(u,v)$-plane.
As $\e$ tends to zero on one end of this moduli
space, the solution tends to the round spherical metric.
As $\e$ tends to zero at the other end, the solution
tends to zero.

\head Spectral Theory for the Delaunay Solutions \endhead
Although the Alexandrov reflection argument shows that we have already
described the full moduli space $\Cal M_{\Lambda}$ when $\Lambda$
has only two elements, we proceed further here to analyze the
linearization of the scalar curvature operator around a Delaunay
solution $g_{\e}$. This case will serve as the model,
and an important ingredient, for the more general linear analysis
later.  Thus, for $v$ a suitably small function, let
$$
\aligned N_{\e}(v) & = \Lap_\e(1+v) - \frac{n(n-2)}4 (1+v)
+ \frac{n(n-2)}4 (1+v)^{\frac {n+2}{n-2}} \\
& = \Lap_\e v + n v + Q(v) \endaligned \tag 2.11
$$
where $Q(v)$ is the nonlinear, quadratically vanishing term
$$
Q(v) = \frac {n(n-2)}4\left( (1+v)^{\frac{n+2}{n-2}} - 1
- \frac{n+2}{n-2}v\right). \tag 2.12
$$
Solutions of $N_{\e}(v) = 0$ correspond to other complete metrics
on the cylinder with scalar curvature $n(n-1)$, and hence correspond to
other Delaunay solutions. The linearization of $N_\e$ about $v=0$ is
thus given by
$$
\left.\frac{d\,}{d\sigma}N_{\e}(\sigma \phi)\right|_{\sigma=0} =
L_{\e} \phi = \Lap_\e \phi + n\phi. \tag 2.13
$$
Frequently we shall omit the $\e$-subscript from $L_\e$ when the context
is clear.

Using (2.6) and (2.13) and introducing an eigendecomposition
$\{\psi_j, \lambda_j\}$ for $\Lap_\th$ on $\SS^{n-1}$ we decompose $L$
into a direct sum of ordinary differential operators
$L_j$ with periodic coefficients on $\RR$. The spectral analysis of each
$L_j$, hence of $L$ itself, is accomplished using Bloch wave theory,
cf\. [RS]. One conclusion is that $\text{spec}(L)$ is
purely absolutely continuous, with no singular continuous or
point spectrum. Moreover, each $L_j$ has spectrum arranged into bands,
typically separated by gaps; $L$ itself has spectrum which is
the union of all of these band structures:
$$
\text{spec}(L) = \bigcup_j \text{spec}(L_j).
$$
{}From (2.13) it is clear that $\text{spec}(-L)$ is bounded below by $-n$.
We proceed to analyze this spectrum near $0$.

The first question is to understand the Jacobi fields, i\.e\. elements
of the nullspace of $L$. Any solution of $L\psi = 0$ may be
decomposed into its eigencomponents
$$\psi = \sum a_j(t)\,\psi_j(\theta),$$
where the $a_j$ solve $L_j(a_j) = 0$. Some of these solutions,
although not necessarily all, may be obtained as derivatives of
one-parameter families of solutions of $N(v)=0$. It is common in
geometric problems, cf\. [KKS], for the Jacobi fields
corresponding to low values of $j$ to have explicit geometric
interpretations as derivatives of special families of solutions of the
nonlinear equation. This is the case here, at least for the
eigenvalues $\lambda_0 = 0$ and $\lambda_1 = \dots = \lambda_{n} = n-1$.

For $j=0$ we look for families of solutions of $N(v) = 0$
which are independent of the $\SS^{n-1}$ factor.
There are two obvious examples, one corresponding to
infinitesimal translations in $t$ and the other corresponding to
infinitesimal change of Delaunay parameter:
$$
\aligned  \eta \rightarrow \frac{u_{\e}(t + \eta)}{u_{\e}(t)}
 & \equiv \Phi_1(t,\e;\eta)
\ \text{and} \\ \eta \rightarrow \frac{u_{\e + \eta}(t)}{u_{\e}(t)}
&\equiv \Phi_2(t,\e;\eta).
\endaligned \tag 2.14
$$
Now define
$$
\phi_j = \phi_{j,\e} = \frac {d}{d\eta} \left. \Phi_j(t,\e;\eta)\right|
_{\eta=0},\quad j = 1,2,  \tag 2.15
$$
so that $L_\e \phi_j = 0,\ j = 1,2$. Of course, neither $\phi_1$ nor
$\phi_2$ are in $L^2$ since $L$ has no point spectrum.
In fact, differentiating
$ u_\e(t+T(\e)) = u_\e(t)$, first with respect to $t$ and then
with respect to $\e$, shows that
$$
\phi_1(t+T(\e)) = \phi_1(t), \qquad \phi_2(t+T(\e)) + \phi_1(t) T'(\e)
= \phi_2(t). \tag 2.16
$$
The first of these equalities states that $\phi_1$ is periodic while
the second shows that $\phi_2$ increases linearly, at least so long as
$T'(\e) \ne 0$. Hence the $\phi_j$ are generalized eigenfunctions
for $L$ with eigenvalue $\lambda = 0$; their existence and
slow growth imply that $0$ is in the essential spectrum of $L_0$,
hence that of $L$.
It is also not hard to see that $\phi_1$ and $\phi_2$ are linearly
independent when $\e \ne \bar u$. In fact, if we translate $u_{\e}$ so
that it attains a local maximum at $t=0$, then $\phi_1(0) = 0$,
whereas $\phi_2(0) = (d\e/d\e)/\e = 1/\e$, so that these functions
are not multiples of one another. When $\e= \bar u$, $\phi_1 \equiv 0$,
but $\phi_2$ does not vanish since $\phi_2(0) = 1/\bar u$.
A second solution of $L\phi = 0$
in this case is obtained by translating $\phi_2$ by any
noninteger multiple of its period.  We shall prove below that these
two functions form a basis, for each $\e$, for all temperate solutions
of $L\phi = 0$ on the cylinder; any other solution of this
equation must grow exponentially in one direction or the
other.

To find Jacobi fields corresponding to the first nonzero eigenvalue
of $\Lap_{\th}$, it is easiest first to transform the background cylindrical
metric to the round spherical metric.
The conformal invariance of the equation means that
one-parameter families of solutions may be obtained from
one-parameter families of conformal maps of $\Sn$.
For example, $\Phi_1$ in (2.14) already corresponds to the composition of
$g_{\e}$ with the family of conformal dilations fixing
the two singular points.
We can also consider composition with a family of parabolic
conformal transformations which fixes one or the other of the singular
points. Derivatives of these families lead to the
new Jacobi fields. To compute them we first
transform the equation yet again so that the background
metric is the flat Euclidean metric on $\RR^n$ and the singular
points are at $0$ and $\infty$. The equation becomes
$$
\Lap_{\RR^n} w + \frac{n(n-2)}{4} w^{\frac{n+2}{n-2}} = 0, \tag 2.17
$$
where the function $w$ on $\RR^n$ is related to $u$ on the cylinder
by the transformation
$$
v(t,\th) = e^{(n-2)t/2}\,w(e^{-t}), \quad w(\rho,\th) = \rho^{(2-n)/2}\,
v(-\log \rho,\th),  \tag 2.18
$$
in terms of the polar coordinates $(\rho,\th)$ on $\RR^n$. The Delaunay
solutions correspond to functions $w_{\e} = w_{\e}(\rho)$. A parabolic
transformation of $\Sn$ fixing one of the singular points corresponds
to a translation of $\RR^n$ if that singular point corresponds to
the point at infinity. Infinitesimal translation is just differentiation
by one of the Euclidean coordinates $x_j$. But
$$
\frac{\del w_{\e}}{\del x_j} = \th_j \frac{\del w_{\e}}{\del \rho},
$$
and $\th_j = x_j / \rho$ is an eigenfunction for $\Lap_{\th}$
with eigenvalue $n-1$. So $\del w_{\e}(\rho)/\del \rho$ is a solution
of $L_1 \phi = 0$, where $L_1$ is written in terms of $\rho$ instead
of $t$. Transforming back to the cylinder once again we arrive
at the solution
$$
\aligned
\phi_3(t) &=  e^t \left( \frac{u'(t)}{u(t)} + \frac{n-2}2 \right) \\
 & = e^t \left( \phi_1(t) + \frac{n-2}2\right).
\endaligned
\tag 2.19
$$
To obtain the other, exponentially decreasing, solution to $L_1 \psi = 0$,
we simply observe that inversion about the unit sphere
in $\RR^n$ corresponds to reflection in $t$ about $t=0$.
If $u_{\e}$ is positioned so that it assumes a maximum
at $t=0$, then $u_{\e}(-t) = u_{\e}(t)$, and we obtain from (2.19)
the solution
$$
\phi_4(t) = e^{-t} \left( \frac{n-2}2 - \phi_1(t) \right). \tag 2.20
$$
The solutions $\phi_1 \dots \phi_4$ are also recorded
and employed in [AKS].

It is also important to have some understanding of the continuous
spectrum near $0$; we do this for each $L_j$ individually.
\proclaim{Lemma 2.21} For every $\e \in (0, \bar u]$ and
$j \ge 1$, $0 \notin \text{spec}(L_j)$. Hence the spectral analysis
of $L$ near the eigenvalue $0$ reduces to that of $L_0$. The
only temperate solutions of $L_0 \phi = 0$ are precisely the
linear combinations of the solutions $\phi_1$ and $\phi_2$ of
(2.15).
\endproclaim
\demo{Proof} The eigenvalues of $\Lap_\th$ are $\lambda_0 = 0,\
\lambda_1 = 1-n$, etc. For $j > n$ (numbering the eigenvalues
with multiplicity, as usual) $\lambda_j < -n - 1$. Thus
$$
\aligned L_j &= u^{-\frac{4}{n-2}}\del_t^2 + 2u^{-\frac{n+2}{n-2}}
u_t \del_t + u^{-\frac{4}{n-2}}\lambda_j + n \\
& < u^{-\frac{4}{n-2}}\del_t^2 + 2u^{-\frac{n+2}{n-2}}
u_t \del_t + n(1 - u^{-\frac{4}{n-2}}) -
u^{-\frac{4}{n-2}} \endaligned \tag 2.22
$$
Since $0 < u \le 1$ for every $\e$ and $t$, the term of order
zero in this expression is always strictly negative. Standard ODE
comparison theory now implies that an arbitrary solution of $L_j \phi = 0$
must grow exponentially either as $t \rightarrow +\infty$ or $-\infty$.
Clearly then $0 \notin \text{spec}(L_j)$ for $j > 1$.

The same conclusion is somewhat more subtle when $j = 1$.
First of all, observe that because $\phi_3$ and $\phi_4$ above
span all solutions to the ODE $L_1 \psi = 0$, and since both of
these functions grow exponentially in one direction or the
other, it is evident that $0$ is not in the spectrum of $L_1$.
We wish to know, though, that $0$ is actually below all
of the spectrum of $-L_1$. This is no longer obvious since
the term of order zero in
$$
L_1 = u^{-\frac{4}{n-2}}\del_t^2 + 2u^{-\frac{n+2}{n-2}}
u_t\del_t + (1-n)u^{-\frac{4}{n-2}} + n \tag 2.23
$$
is no longer strictly negative. Conjugate
$L_1$ by the function $u^p$, where $p$ is to be chosen. A simple
calculation gives
$$
u^{-p}L_1 u^p = u^{-\frac{4}{n-2}}\del_t^2 + (2p + 2)uu_t\del_t + A,
$$
where
$$
A = u^{-\frac{2n}{n-2}}\left\{ p(p+1)u_t^2 +
\left(\frac{p(n-2)^2}4 + 1-n\right) u^2
+ \left( n - \frac{pn(n-2)}4 \right) u^{\frac{2n}{n-2}} \right\}.
$$
Somewhat miraculously, upon setting $p = 2/(n-2)$ this expression
reduces to
$$
\gathered A = \frac{2n}{(n-2)^2} u^{-\frac{2n}{n-2}}
\left( u_t^2 - \frac{(n-2)^2}4 u^2 + \frac{(n-2)^2}4 u^{\frac{2n}{n-2}}
\right) \\
= \frac{4n}{(n-2)^2} u^{-\frac{2n}{n-2}} H(\e), \endgathered
$$
where $H(\e)$ is the Hamiltonian energy (2.3) of the solution $u_{\e}(t)$.
This is a negative constant for every $\e \in (0, \bar u]$,
and so once again we have obtained an ordinary differential operator
with strictly negative term of order zero. $u^{-2/(n-2)} L_1
u^{2/(n-2)}$ is unitarily equivalent to $L_1$, and by the previous
argument, the spectrum of the first of these operators is strictly
negative, and in particular does not contain $0$. Hence the
same is true for $L_1$ and the proof is complete.
\enddemo

We can say slightly more about the spectra of the $L_j$. In fact,
when $\e = \bar u$,
$$ L_j = \frac {n}{n-2}\left(\del_t^2 + \lambda_j \right) + n, $$
and
$$ \text{spec}(-L_j) = \left[-n\left(1 + \frac{\lambda_j}{n-2}\right),
\infty\right).$$
In particular, $\text{spec}(-L_0) = [-n, \infty)$ and $\text{spec}(-L_1)
= [n/(n-2), \infty)$.  The spectrum of $-L_{\bar u}$ is the union of
these infinite rays. As $\e$ decreases, gaps appear in these
rays, and bands start to form (in fact, by (2.8) the first gap appears
at $-3n/4$). As $\e$ decreases, the first band extends from $-n$ to
some value to the left of $-3n/4$; the second band begins somewhere to
the right of this point, and always ends at $0$. It is not important
for our later work that $0$ is always on the end of the second
band (or at least, at a point where the second Bloch band function
has a turning point), but it is rather amusing that we may determine
this. The explanation is that from (2.16), whenever $T'(\e) \ne 0$ (which
certainly happens for almost every $\e < \bar u$),
the Jacobi field $\phi_2$ grows linearly. This signifies
that $0$ is at the end of a band, or at the very least, at a
turning point for one of the Bloch band functions for $L_0$ which
parametrize the bands of continuous spectrum.
As $\e$ continues to decrease to $0$ these bands shrink further and
further. It is more revealing to consider the operator $L_\e$ written
in the geodesic coordinates $(r,\th)$ as in (2.7). As noted earlier, as $\e
\rightarrow 0$ this operator remains periodic, but develops
singularities every $\pi$ units. Geometrically, the metrics $g_\e$
are converging to an infinite bead of spheres of fixed scalar
curvature $n(n-1)$. The spectrum of $-L$ for this limiting metric
is a countable union of the spectrum of $-\Lap - n$ on one
of these spheres, i\.e\. consists of a countable number of isolated
points $\{-n, 0, n+2, \dots\}$, each with inThe bands of continuous
spectrum have coalesced into these
infinite-multiplicity eigenvalues. Note that in this limit,
the infimum  of $\text{spec}(L_1)$ has increased to zero at
a rate which may be estimated by a power of $\e$.

\specialhead III. Generalities on Solutions with Isolated Singularities
\endspecialhead

In this section we collect various results concerning solutions
of the singular Yamabe problem with isolated singularities, particularly
in the context of conformally flat metrics, which will be required later.

\head Asymptotics \endhead
The Delaunay solutions discussed in the last section are interesting
explicit solutions in their own right. However, their importance
is due to the fact that they are asymptotic models for arbitrary
solutions of the singular Yamabe problem, at least in the
conformally flat setting. Before stating this result more carefully we
introduce some notation. The half-cylinder $(0,\infty)_t \times
\SS^{n-1}_{\th}$ is conformally equivalent to the punctured ball
$\Bbb B^n \backslash \{0\}$; an explicit conformal map is given
by sending $(t,\theta)$ to the polar coordinates $(\rho,\theta)$,
where $\rho = e^{-t}$. A Delaunay solution $u_{\e}(t)$ may be
transformed to the function
$$ \tilde u_{\e}(\rho) = \rho^{\frac{2-n}2} u_{\e}(-\log \rho). \tag 3.1 $$
Using the conformal invariance of the conformal Laplacian, this function
solves
$$ \Lap \tilde u + \frac{n(n-2)}4 \tilde u^{\frac{n+2}{n-2}} = 0.
\tag 3.2 $$

\proclaim{Theorem 3.3} Let $\tilde u \in \Cal C^{\infty}(\Bbb B^n
\backslash \{0\})$ be a positive solution of equation (3.2) with a
nonremovable singularity at the origin. Then there exists
some $\e \in (0, \bar u]$ and $A>0$ so that $\tilde u$ is asymptotic to
the modified Delaunay solution $\tilde u_{\e}$ in the sense that
$$ \tilde u(\rho,\theta) = (1 + O(\rho^{\a}))\tilde u_{\e}(A\rho, \theta),
\tag 3.4 $$
for some $\a > 0$; $A$ corresponds to a translation of the $t$
parameter. The analogous estimate still holds whenever $\tilde u$
and $\tilde u_{\e}$ are replaced by any of their derivatives
$(\rho \del_{\rho})^j\del_{\th}^{\b}$.
\endproclaim
This result was proved by Caffarelli, Gidas and Spruck [CGS]
using a fairly general and complicated form of the Alexandrov reflection
argument. They actually only give a somewhat weaker estimate where
$O(\rho^{\a})$ is replaced by $o(1)$. An alternate
argument, relying on more direct geometric and barrier methods,
was obtained at around the same time by Aviles, Korevaar and Schoen
in the unpublished work [AKS], and they obtained the stronger form.
This last work is very close in method to the proof of the analogous
result for complete constant mean curvature surfaces embedded in $\RR^3$
presented in [KKS].
It is possible to calculate the decay rate $\a$ in terms of spectral
data of the operator $L_{\e}$, and we shall indicate this argument
in the next section. To our knowledge it is unknown whether some form
of this result has been proved when the background metric is not
conformally flat, but it seems likely to be true.

Note that
the estimate (3.4) may be restated for the transformed function
$u(t,\th) = e^{-(n-2)t/2} \tilde u(e^{-t}, \th)$ on the cylinder. Now
$$ u(t,\th) = (1 + O(e^{-\a t})) u_\e (t). \tag 3.5 $$
This is the form we shall use. It states that an arbitrary solution
on the half-cylinder converges exponentially to a Delaunay solution.

\head Poho\v zaev Invariants \endhead
A key ingredient in Schoen's construction [S2] of solutions of
(1.3) with isolated singularities is his use of balancing conditions
for approximate solutions. These conditions follow from the
general Poho\v{z}aev identity proved in [S2]. In the present
setting, when $(\Omega,g)$ is a compact, conformally flat manifold
with boundary, $T$ is the trace-free Ricci tensor of $g$,
$X$ is a conformal Killing field on $\Omega$,
and $R(g)$ is constant, then
$$
\int_{\del \Omega} T(X, \nu) d\sigma = 0.  \tag 3.6
$$
Here $\nu$ is the outward unit normal to $\del \Omega$ and $d\sigma$
is surface measure along this boundary. There is a more general
formula involving an integral over the interior of $\Omega$ when
the scalar curvature is not constant, or when $X$ is not conformal
Killing.

When $\Omega \Subset \SS^n \backslash \{p_1, \dots , p_k\}$, there
are many conformal Killing fields to use.  An important class of
these are the `centered dilations.'  Any such $X$
is equal to the gradient of the restriction of a linear function
$ \ell(q) = \langle q, v \rangle$, where $v \in \RR^{n+1}$ and
$\SS^n \subset \RR^{n+1}$ is the standard embedding. Explicitly,
$$ X_q = v - \langle q, v \rangle q,\qquad q \in \SS^n, \tag 3.7 $$
since the unit normal to $T_q\SS^n$ is just $q$. Let $\frak c$ denote
the set of all such $X$. It is a subspace of the Lie algebra,
$\frak o (n+1,1)$, of the conformal group $O(n+1,1)$.

By (3.6), each $X \in \frak  o (n+1,1)$ determines
an element of $H^{n-1}(\SS^n \backslash \Lambda, \RR)$, associating
to a hypersurface $\Sigma$ the number
$$
\Cal P(\Sigma, X;g) = \int_{\Sigma} T(X,\nu)\,d\sigma \tag 3.8
$$
The dual homology space is generated by the classes of hyperspheres
$\Sigma_i$, where, for each $p_i \in \Lambda$,
$\Sigma_i=\partial B_{r}(p_{i})$ for $r$ sufficiently small, is chosen
so that no other $p_j$ is in the
same component of $\SS^n \backslash
\Sigma_i$ as $p_i$. These classes satisfy the single relation
$[\Sigma_1] + \cdots + [\Sigma_k] = 0$. The number in (3.8) will be written
$\Cal P_i(X;g)$ (or simply $\Cal P_i(X)$) when $\Sigma = \Sigma_i$.
Alternately, as suggested
by this notation, we shall also regard each $\Cal P_i$ as a linear
functional on the $X's$, i\.e\. as an element of $\frak o^*(n+1,1)$.
These functionals will be called the Poho\v{z}aev invariants of
the solution metric $g$. They satisfy
$$ \Cal P_1 + \cdots + \Cal P_k = 0.  \tag 3.9 $$

Theorem 3.3 makes it possible to calculate at least some components
of the $\Cal P_i$, by computing the Poho\v{z}aev invariants for
the Delaunay solutions; we do this now.
Since $H^{n-1}(\SS^n \backslash \{p_1, -p_1\},\RR)
\equiv \RR$, it suffices to compute the single number $\Cal P_1(X)$
for each $X \in \frak c$. As described above, $X$ may be identified
with an element $v \in \RR^{n+1}$. A straightforward calculation,
similar to one given in [P2],
shows that this invariant is, up to a constant, simply
the Hamiltonian energy $H(\e)$ of the particular Delaunay solution:
\proclaim{Lemma 3.10} Using the identification of $X \in \frak c$ with
$v \in \RR^{n+1}$, the Poho\v{z}aev invariant for the Delaunay
metric $g_{\e}$ on $\SS^n \backslash \{p_1, -p_1\}$ is given by
$$ \Cal P_1(X) = c_n\,H(\e)\,\langle v, p_1 \rangle. $$
Here $c_n$ is a nonvanishing dimensional constant (identified
explicitly in [P2]).
\endproclaim

The asymptotics in Theorem 3.3 (even with the sharp estimate
of $\a$ we will give later) do not give enough information for the
invariant $\Cal P_i$ to be computed for every $X$; the point is that
if $X$ has associated vector $v$ perpendicular to $p_i$ then Lemma 3.10
shows that there is no `formal' asymptotic contribution to the invariant,
but unfortunately the decay is not sufficient for there to be
no `perturbation' contribution. The one case where this is not an issue
is when $X$ is chosen to have associated vector exactly $p_i$. The
associated invariant $\Cal P_i(X)$ in this case will be called
the {\it dilational} Poho\v{z}aev invariant and denoted $\Cal D_i$
(or $\Cal D_i(g)$)
\proclaim{Corollary 3.11 ([P2])} The dilational Poho\v{z}aev invariant
$\Cal D_i(g)$ associated to the puncture $p_i$ and the solution metric
$g$ is equal to $c_n\,H(\e_i)$, where $\e_i$ is the Delaunay parameter
giving the asymptotic model for the metric $g$ at $p_i$, as provided
by Theorem 3.3.
\endproclaim

\head Local Compactness Properties of the Moduli Space \endhead
It is natural to determine the possible ways that solution metrics
$g$ of our problem can degenerate. Phrased more geometrically,
we wish to determine what the ends of the
moduli space $\Cal M_{\Lambda}$ look like, and to determine
a geometrically natural compactification. Any $\Cal M_{\Lambda}$,
when $\Lambda$ contains just two elements, may be identified
with any other, and we may call this space simply $\Cal M_2$.
In this space, all $g \in \Cal M_2$ are Delaunay, and indeed
$\Cal M_2$ may be identified with the open set $\Omega \subset
\RR^2$, $\Omega = \{H < 0\}$ as described in \S2. This space is two
dimensional, smooth, and locally compact. A sequence of
elements in $\Cal M_2$ degenerate only when the
Delaunay parameter $\e$ (or `neck-size' of the solutions) tends
to zero. Recalling that the Hamiltonian energy $H$ is strictly monotone
decreasing for $\e \in (0, \bar u]$, this may be restated as saying
that Delaunay solutions degenerate only when $H(\e) \rightarrow 0$.

A similar statement is true when $\Lambda$ has more than two
elements. As we have seen in Corollary 3.11, the
Poho\v{z}aev invariants $\Cal P_i$, and indeed just the dilational
Poho\v{z}aev invariants $\Cal D_i$, determine the Hamiltonian
energies of the model Delaunay solutions for the metric $g$ at
each puncture $p_i$. Clearly, then, if there exists some
particular $p_i$ such that for a sequence of solution
metrics $g_j \in\Cal M_{\Lambda}$ the dilational Poho\v{z}aev
invariant $\Cal D_i(g_j)$ tends to zero as $j \rightarrow \infty$,
then this sequence should be regarded as divergent
in $\Cal M_{\Lambda}$. The complementary statement, generalizing
the situation for $\Cal M_2$, is also still valid, and was
proved by Pollack [P2]:

\proclaim{Proposition 3.12} Let $g_j$ be a sequence of metrics
in $\Cal M_{\Lambda}$, such that for each $i=1,\ldots,k$ the
dilational Poho\v{z}aev invariants
$\Cal D_i(g_j)$ are bounded away from zero.  Then there is a
subsequence of the $g_j$ converging to a metric $\bar g \in
\ML$.  The convergence is uniform in the $\cin$ topology relative to
$\bar g$, or indeed relative to any of the $g_j$, on compact subsets
of $\SL$.
\endproclaim

The question of what happens to divergent sequences in $\ML$
may also be determined. Before describing this we return to
$\Cal M_2$. There is an `obvious' compactification,
$\overline{\Cal M}_2$, which is the closure of $\Omega$ in $\RR^2$.
$\overline{\Cal M}_2 \backslash \Cal M_2$ decomposes into
two disjoint sets: the first contains the single point $\{0,0\}
\in \RR^2$, which corresponds to the conformal factor $0$,
while the second is the orbit $\{u_0(t), v_0(t)\}$ passing
through $\{1,0\}$, which corresponds to the incomplete spherical
metric on $\Sn \backslash \{p_1, -p_1\}$. These latter points
on the compactification of $\Cal M_2$ may be identified with
the nonsingular round metric on $\Sn$ itself, and these are
themselves metrics of constant scalar curvature $n(n-1)$ on
$\Sn$. $\overline{{\Cal M}_2}$ is a stratified space; its principal
stratum is $\Cal M_2$ itself, the codimension one (and
one dimensional) stratum consists of copies of the nonsingular round
metric on $\Sn$, and finally the codimension two (zero dimensional)
stratum consists of the single trivial solution of the PDE, which
is a completely degenerate metric.

There is a similar compactification of $\ML$ and corresponding
decomposition of $\CML$ when the cardinality of
$\Lambda$
is larger than two (and still finite). The following is a corollary
of the proof of Proposition 3.12:

\proclaim{Corollary 3.13} Let $g_j \in \ML$ be a sequence
such that for some nonempty subset of points
$p_{i_1},\dots, p_{i_s}$, the invariants $\Cal D_{i_1}(g_j),
\dots, \Cal D_{i_s}(g_j)$ are uniformly bounded above by
some $-\eta_0 < 0$ for all $j$, and all other $\Cal D_i(g_j)$
tend to zero. Then this sequence has a convergent
subsequence converging to a metric
${\bar {g}}\in{\Cal{M}}_{\Lambda^{\prime}}$, i\.e\. ${\bar{g}}$  is still
singular at the points $\Lambda' \equiv \{p_{j_1}, \dots, p_{j_s}\}$, but
extends smoothly across the points in $\Lambda \backslash \Lambda'$ and
$R({\bar{g}})=n(n-1)$. If all ${\Cal{D}}(g_{i}, p_{j})$ tend to zero as
$i \rightarrow \infty$ then either $g_i$ tends to zero, uniformly on
compact subsets of $S^{n}\setminus\Lambda$, or else $g_i$ converges to
the round metric on $S^{n}$.
\endproclaim

This means that the compactification
$\CML$ contains copies of $\Cal M_{\Lambda'}$ for
certain subsets $\Lambda' \subset \Lambda$. In addition it may
contain  a set
whose points correspond to copies of the round metric
on $\Sn$, and finally a set whose points correspond to
the `zero metric'. Later in this paper we shall give
a somewhat better description of $\CML$ once we
have determined the structure of $\ML$ itself better.
Notice also that it must be somewhat subtle to determine
precisely which subsets $\Lambda' \subset \Lambda$ have
$\Cal M_{\Lambda'}$ occurring in $\overline{\ML}$.
The subsets $\Lambda'$ which can arise in the description above
are determined by the Poho\v{z}aev balancing condition (3.9).
For example,
when $\Lambda$ has two elements, there is no piece of the
boundary corresponding to $\Cal M_{\Lambda'}$ for $\Lambda'$
having just one element; this is because no complete solution
of our problem exists on $\Sn \backslash \{p\}$. {}From the present
perspective, this holds because $\Cal D_1(g_j) -
\Cal D_2(g_j) = 0$ for any $j$ (here $\Lambda = \{p_1,p_2\}$)
by virtue of (3.9), so that if one of these numbers tends to
zero, the other must as well.  (The difference in the signs
here from (3.9) is because in (3.9) the same conformal
Killing field $X$ is used in each Poho\v{z}aev invariant,
whereas here we use $X$ corresponding to $p_1$ for one
and $-X$ corresponding to $p_2 = -p_1$ for the other.)
The subtle point is that the complete Poho\v zaev invariants,
rather than just the dilational ones,
are required when $\Lambda$ has more than two points,
and these are not determined just linear algebraically by the
location of the $p_i$ and the Delaunay models at these punctures.

\specialhead IV. Linear Analysis on Manifolds
with Asymptotically Periodic Ends \endspecialhead

In this section we prove various results concerning the analysis
of the Laplacian and the linearization of the scalar curvature
operator about solution metrics $g \in \Cal M_{\Lambda}$. Many of
the basic results hold more generally, e\.g\. for the Laplacian
on manifolds with asymptotically periodic ends.

\head The Linearized Operator \endhead

For the remainder of this paper, $L$ will always denote the
linearization about the constant solution $v=0$ of the nonlinear operator
$$ \aligned
N(v) & = \Lap_g(1+v) - \frac{n(n-2)}4 (1+v) + \frac{n(n-2)}4(1+v)^{
\frac{n+2}{n-2}} \\
& = \Lap_g v + n v + Q(v), \endaligned \tag 4.1 $$
where $Q(v)$ is the same as in (2.12), so that, as before,
$$ L = \Lap_g + n. \tag 4.2 $$

Let $E_j$ denote a neighborhood of the puncture $p_j$ which
is conformally equivalent to a half-cylinder $[0,\infty)_t \times \SS^{n-1}
_\th$. We fix these cylindrical coordinates around each $p_j$.
By virtue of Theorem 3.3, $L$ can be treated on each $E_j$ as
an exponentially small perturbation of the corresponding operator
$L_{\e_j}$ for the periodic Delaunay metric $g_{\e_j}$ which is the
asymptotic model for $g$ on $E_j$.  Thus on each $E_j$ we may write
$$ L = L_{\e_j} + e^{-\a t}F,  \tag 4.3 $$
where $F$ is a second order operator with coefficients bounded
in $\cin$ as functions of $(t,\th)$.

The linear analysis of the Laplacian on manifolds with
asymptotically periodic ends is remarkably similar to that
for manifolds with asymptotically cylindrical
ends, as detailed for example in [Me].  In particular, the Fredholm
theory for such an operator on exponentially weighted Sobolev (or H\"older
or ...) spaces has an almost identical statement in either
case, although the proofs are rather different. The Fredholm theory
in this asymptotically periodic setting was previously developed
by Taubes in [T], although we proceed somewhat further into the
linear analysis as we need more detailed results. It is also possible
to develop a full scattering theory for the Laplacian on
these manifolds. We shall actually require and develop some
scattering theoretic results to clarify the nature of the moduli space
$\Cal M_{\Lambda}$.

\head The Fourier-Laplace Transform \endhead

The basic tool for the parametrix construction, upon which all the
linear analysis relies, is the Fourier-Laplace transform, in a
form employed already by Taubes [T]. We proceed to develop some properties
of this transform, here denoted $\Cal F$.

The function spaces we shall use here are exponentially
weighted Sobolev spaces based on $L^2(\SL;dV_g)$; these will be written
$H^s_{\ga}(\SS^n \backslash \Lambda)$, or just $H^s_{\ga}$, for
$\ga, s \in \RR$, $s > n/2$. The last condition ensures
that the spaces behave well under nonlinear operations. To define
them, decompose $\SS^n \backslash \Lambda$ into the union of the
ends $E_1, \dots, E_k$ and a compact piece $K$. Over $K$ an element
$h \in H^s_{\ga}$ restricts to an ordinary $H^s$ function.
Over $E_j$, $h = e^{\ga t}\tilde h$, where $\tilde h \in H^s([0,\infty)
\times \SS^{n-1}, dt\,d\th)$. Note that the measure here is uniformly
equivalent to the one induced by $g$ for any $g \in \ML$.

The transform $\Cal F$ is, strictly speaking, defined for functions
on the whole cylinder $C = \RR \times \SS^{n-1}$. It would be somewhat
more natural to first develop its properties acting on e\.g\. the
Schwartz space $\Cal S$, but we shall specialize immediately
to functions with support on just one half of $C$.
So let $h(t,\th)\in H^s_\ga$ on $E_j$, and assume $h = 0$ for $t \le 1$.
Set
$$ \hat h(t,\z,\theta) =\Cal F(h) = \sum_{k = -\infty}^{\infty} e^{-ik\z}
h(t + k,\theta).  \tag 4.4 $$
Assume for the moment that $h$ is smooth. Then, since $h$ decays like
$e^{\ga t}$, this series converges provided $\Im \z \equiv \nu <
-\gamma$. We have set $\z = \mu + i\nu$, so that $\Re \z = \mu$.
$\hat h(t,\z,\theta)$ depends holomorphically
on $\z$ in the region $\nu < -\ga$. When $h$ is only assumed
to be in $H^s_{\gamma}$, $\hat h(t,\z,\theta)$ will still depend
holomorphically on $\z$ in the same region, but as a function with
values in the space $H^s$. $\hat h$ is continuous in $\nu \le -\ga$ as
a function with values in $H^s$. These results follow from the Plancherel
formalism.

The tranform $\Cal F$ is invertible, and its inverse is
given by contour integration. To make the following equation
clearer, assume that $t \in \Bbb R$, and $\tilde t$ is its reduction
mod 1 (so that $0 \le \tilde t < 1$). Then when $\ell \le t < \ell + 1$
so that we may write $t = \tilde t + \ell$,
$$ h(t,\theta) = \frac 1{2\pi} \int_{\mu = 0}^{2\pi} e^{i \ell \z}
\hat h(\tilde t,\z,\theta) d\mu.  \tag 4.5 $$
In this formula we integrate along a line $\Im \z = \nu_0$.
By Cauchy's theorem this contour may be shifted to allow $\nu_0$ to
be any number less than $-\gamma$. If, as we are assuming, $h$ vanishes
for $t < 0$, then $\hat h(\tilde t,\z,\theta)$ not only extends
to the lower half $\z$ plane, but decays like $e^{\nu}$ there.
Shifting the contour arbitrarily far down shows that the integral
(4.5) vanishes for any $\ell < 0$, as it should.
By a similar argument, if $h(t,\theta)$ is defined by the integral (4.5)
taken along some contour $\Im \z = \nu_0$, where the integrand
$\hat h(\tilde t, \z, \theta)$ is only assumed to be defined
along that line, then $h \in H^s_{-\nu_0}$. In particular, if $\hat h$
is holomorphic in some lower half plane $\nu < -\ga$, and continuous
with values in $H^s$ as a function of $\tilde t$ up to this upper
boundary, then $h \in H^s_{\ga}$.

Next, reindexing the sum defining $\hat h$ gives
$$ \hat h(t+1, \z,\theta) = e^{i\z} \hat h(t,\z, \theta). \tag 4.6 $$
This just means that $\hat h(t,\z,\theta)$ is a section of the flat bundle
on $\SS^1 \times \SS^{n-1}$ with holonomy $\z$ around the $\SS^1$
loop. This bundle is isomorphic to the flat bundle with
trivial holonomy; the bundle map is given by conjugating
by $e^{i\z t}$. Thus the function
$$ \tilde h \equiv e^{-i\z t} \hat h(t,\z,\theta) e^{i\z t}
\tag 4.7 $$
satisfies $\tilde h(t+1,\z,\theta) = \tilde h(t,\z,\theta)$.

\head Fredholm Theory \endhead

The basic Fredholm result for the linearization $L$ may
now be stated and proved.
\proclaim{Proposition 4.8} There exists a discrete set of numbers
$\Gamma \subset \RR$ such that the bounded operator
$$
L:  H^{s+2}_{\gamma}(M) \longrightarrow H^{s}_{\gamma}(M) \tag 4.9
$$
is Fredholm for all values of the weight parameter $\gamma \notin\Gamma$.
In particular, $0 \in \Gamma$, so the map (4.9) is not Fredholm when
$\gamma = 0$, i\.e\. on the ordinary unweighted Sobolev spaces, but
is Fredholm for all values of $\gamma$ sufficiently near, but not
equal to zero.
\endproclaim
\demo{Proof}  It suffices to construct a parametrix $G$
for $L$ so that $LG - I$ and $GL - I$ are compact operators
on $H^s_{\ga}$ and $H^{s+2}_{\ga}$, respectively.  As usual,
$G$ may be constructed separately on each piece of the decomposition
of $\SL$ into a compact piece and the $k$ ends $E_1, \dots E_k$
around each $p_j$. The parametrix construction on the compact piece
is the standard microlocal one since $L$ is elliptic. We construct
a parametrix on each $E_j$ using the Fourier-Laplace transform.

Fix one of these ends, $E_j$, and let $g_\e$ be the model
asymptotic Delaunay metric for the fixed metric $g$ there.
The corresponding model operator $L_{\e}$, as in (4.3),
has periodic coefficients (for notational convenience we assume that
the period is one here) so it acts on sections of the flat
bundle with holonomy $\z$ described above by the obvious rule
$\widehat (L_\e h) \equiv \hat L_{\e} \hat h$. This induced
operator $\hat L_{\e}$ looks just like $L_{\e}$ in local
coordinates $(t,\theta)$. This step is the same as conjugating
$L_\e$ by $\Cal F$. We may proceed further and conjugate
$\hat L_{\e}$ by $e^{i\z t}$ so as to act on the
trivial flat bundle. This final induced operator, which depends
holomorphically on $\z$, will be called $\tilde L_{\e}(\z)$:
$$
\tilde L_\e(\z) = e^{-i\z t} L_\e e^{i\z t}. \tag 4.10
$$

The main point of the proof is that $\tilde L_{\e}(\z)$ has
an inverse (say on $L^2(\SS^1 \times \SS^{n-1})$) which depends
meromorphically on
$\z$. This will be a direct consequence of the analytic
Fredholm theorem, which is proved for example in [RS] and states
that a strongly holomorphic family of Fredholm operators, depending on the
complex variable $\z$, either fails to be invertible for every $\z$, or
else is invertible for all $\z$ except for those in some discrete set
in the parameter space.  To check that this result is applicable, simply
note that $\tilde L_{\e}(\z)$ is elliptic for every $\z$, (as already
observed) depends holomorphically on $\z$ and, by the
various transformations we performed, acts on a {\it fixed} function space
(on a fixed bundle) over a compact manifold. Hence it forms a family
of Fredholm operators; standard considerations show that the
holomorphic dependence of the coefficients of this operator on $\z$
ensure the strong holomorphy of the family.

Our next task is to show that $\tilde L_{\e}(\z)$
is invertible for some value of $\z$. Once this is accomplished,
the analytic Fredholm theorem will imply that the set of poles
$\{\z_j\}$ of $\tilde L_{\e}(\z)^{-1}$ is discrete in
$\Bbb C$. The invertibility of $\tilde L_\e(\z)$ at some  $\z$
is equivalent to the invertibility of $\hat L_\e$ acting
on the flat bundle with holonomy $\z$. When $\z$ is real,
this operator is self-adjoint, hence invertible so long as it
has no nullspace. A solution of $\hat L_\e\psi = 0$
on this bundle lifts to a function $\psi(t,\theta)$ on $\RR \times \SS^{n-1}$
satisfying $\psi(t+1,\theta) = e^{i\z} \psi(t,\th)$.  Since
$\z$ is real here, this lift is bounded and `quasi-periodic.'
However, we have already shown in \S2 that such a function can exist
only when $\z = 0$; in fact, the arguments there show that any temperate
solution must be constant on the cross-section $\SS^{n-1}$, i\.e\.
be a solution of the reduced operator $L_0$. In addition, the functions
$\phi_1$ and $\phi_2$ of (2.15) are independent solutions of this
operator, hence span the space of all solutions since $L_0$ is a
second order ODE. Finally, neither of these transform by a factor $e^{it\z}$
over a period except for $\z = 0$.  Thus $\tilde L_\e$ is invertible
for every real $\z \in (0,2\pi)$. This proves that $\tilde
L_\e(\z)^{-1}$ exists and depends meromorphically on $\z$.

We make some remarks about the set of poles $\{\z_j\}$ of this meromorphic
family of inverses.  First, if $\z$ is in this set, then so is $\z +
2\pi\ell$ for any $\ell \in \Bbb{Z}$. In fact, although
$\tilde L_\e(\z + 2\pi)$ is not equal to
$\tilde L_\e(\z)$, these two operators are unitarily equivalent:
the unitary operator intertwining them is multiplication by the function
$e^{2\pi i t}$. So, for some fixed $\z$, either both
of these operators are invertible or neither of them are.
By discreteness and this translation invariance, the set $\{\z_j\}$
has at most finitely many inequivalent (mod $2\pi$) poles in any
horizontal strip $a < \Im \z < b$. Note in particular that the
only poles on the real line occur at the points in $2\pi\Bbb Z$.
By discreteness again, there are no other poles in some
strip $-\e < \Im\z < \e, \Im\z \ne 0$. We shall usually restrict attention
then only to values of $\z$ in the strip $0 \le \z < 2\pi$ (or
equivalently, in $\Bbb C / 2\pi \Bbb Z$); $\Cal P$ will denote the set of
poles in this strip. Finally, note that the adjoint of $\tilde L_{\e}(\z)$
is precisely $\tilde L_{\e}(\bar \z)$, and so $\Cal P$ is invariant
under conjugation. This means that we can list the elements of
$\Cal P$ as follows:
$$
\Cal P = \Cal P_{\e} = \{ \dots, -\z_2, -\z_1,
\z_0 = 0, \z_1, \z_2, \dots \}.
\tag 4.11
$$
Also define
$$
\Gamma = \Gamma_{\e} = \{\ga_j =
\Im \z_j: \z_j \in \Cal P\} = \{ \dots , -\ga_1,
0, \ga_1, \dots \}. \tag 4.12
$$
The $\ga_j$ form a strictly increasing sequence tending to infinity,
and in particular, $\Gamma$ is a discrete set in $\Bbb R$. We shall
provide an interpretation of the nonzero elements of $\Cal P$ in the next
subsection.

Next we unwind this inverse $\tilde L_{\e}(\z)^{-1} \equiv \tilde
G_{\e}(\z)$ to obtain a parametrix for $L_{\e}$ on $E_j$. The inverse
for $\hat L_{\e}$ acting on the flat bundle with holonomy $\z$ is
given by $\hat G = e^{it\z} \tilde G e^{-it\z}$, for every $\z
\notin \Cal P$. To obtain an inverse for $L_{\e}$ we simply
conjugate $\hat G$ with the Fourier-Laplace transform $\Cal F$.
Recall, though, that this makes sense only if the contour integral in the
definition of $\Cal F^{-1}$ avoids $\Cal P$. So, by conjugating, and
in doing so, taking this integral along a contour $\Im\z = -\ga$, for any
$\ga \notin \Gamma$, we obtain an inverse for $L_{\e}$. Since integrating
along this contour produces a function in $H^s_{\gamma}(\Bbb R
\times \SS^{n-1})$, we have obtained an inverse $G = G_{\e,\gamma}$
for $L_{\e}$ acting on $H^s_{\ga}$ functions on $E_j$,
supported in $t \ge 1$ for all $\ga \notin \Gamma$.
Clearly this $G$ is a parametrix with compact
remainder for $L = L_{\e} + e^{-\a t}F$ acting on $H^s_{\ga}(E_j)$.

As a final step, couple these parametrices we have produced on each end
$E_j$ to the interior parametrix to obtain a global parametrix for
$L$ on all of $\SL$ with compact remainder.
Both right and left parametrices of this type may be obtained this way, and
their existence implies immediately that $L$ is Fredholm on
$H^s_{\ga}$ whenever $\ga \notin \Gamma$.
It is a simple exercise to check that $L$ does not even have closed
range when $\ga \in \Gamma$, in particular, when $\ga = 0$.
\enddemo

It is obviously of interest to determine when $L$ is actually
injective or surjective. In general this is a rather complicated
question, but we make note of the following result.
\proclaim{Corollary 4.13} Suppose that $L$ has no $L^2$ nullspace
(i\.e\. is injective on $H^s_{\gamma}$ for $\gamma \le 0$),
then for all $\d > 0$ sufficiently small (with $\delta < \a$)
$$ \align
&L: H^{s+2}_{\delta} \longrightarrow H^s_{\delta} \quad \text{is
surjective},\\
&L: H^{s+2}_{-\delta} \longrightarrow H^s_{-\delta} \quad \text{is
injective}. \endalign $$
\endproclaim
This first statement follows from the second and duality, since
$L$ is self-adjoint on $L^2$, i\.e\. when $\ga = 0$.

\head Asymptotic Expansions \endhead

The Fourier-Laplace transform may also be used to obtain existence of
asymptotic expansions for solutions to $Lw=0$ on each end $E_j$.
For our purposes in this paper, it will only be important
to know that any such $w$ has a leading term in its expansion,
which decays (or grows) at some specified rate, and an error which
decays at a faster rate. One use of this will be to estimate the
exponent $\a$ appearing in (3.5).
However, the full expansion is not much harder to prove, so
we will do this too.

The starting point is the fact, discussed above, that the Fourier-Laplace
transform of a function $w \in H^s_{\ga}(\RR \times \SS^{n-1})$, supported
in $t \ge 0$, say, is a function $\hat w(t,\theta,\z)$ which is
holomorphic in the half-plane $\Im\z < -\ga$ (and taking values
in $H^s(\SS^1 \times \SS^{n-1})$). For a general function of this type,
the most that can be concluded is that it has a limit on the line
$\Im\z = -\ga$ in the appropriate $L^2$-sense discussed earlier.

When the solution $w$ is transformed to a function
$\bar w(r,\th) = w(-\log r,\th)$
on $\Bbb B^n \backslash \{0\}$, elliptic regularity
shows that $\bar w$ is a conormal distribution with respect
to the origin, i\.e\. has stable regularity when differentiated arbitrarily
often with respect to the vector fields $r\del_r, \del_{\th}$
(this is just the same as $w$ itself having stable regularity
with respect to $\del_t,\del_{\th}$); [M1], for example, contains a
discussion of conormal regularity. For degenerate operators of a type
closely related to $L$, one expects solutions to be polyhomogeneous
conormal, cf\. [M1] or [Me].
Polyhomogeneity is simply the property of having an asymptotic expansion
in increasing (possibly complex) powers of $r$ and integral powers of
$\log r$, with coefficients smooth in $\th$ (as functions of $t$ these
expansions are in analogous powers of $e^t$ and $t$).
Alternately, polyhomogeneous conormal distributions
may be characterized as those with Mellin transforms, already defined and
holomorphic in some lower half-plane in $\Bbb C$, extending
meromorphically to the whole complex plane, with only finitely
many poles in any lower half-plane, all of which are of
finite rank.

$\bar w$ will not be polyhomogeneous, except when $\e = \bar u$ and the
underlying metric is cylindrical. In fact, $\bar w$ (cut off so as to be
supported in $r \le 1$) will have Mellin transform still defined and
holomorphic in a half-plane and extending meromorphically to all of $\Bbb C$,
but now its poles are arranged along lattices on
a countable discrete set of horizontal lines. Each pole is still
finite rank.  This meromorphic structure of the Mellin transform of $\bar w$
is equivalent to the fact that the terms in the expansion for $w$ have
the form $e^{-\b t}w_j$, where $w_j$ is periodic in $t$. The residues at
the poles along a fixed horizontal line at height $\b$ will be the
Fourier series coefficients of $w_j$.

This more general discussion has involved the Mellin transform
as $r \rightarrow 0$ (or equivalently, the standard Fourier transform
as $t \rightarrow \infty$). We shall revert now to the Fourier-Laplace
transform; it is much less flexible than either of
the other transforms, since it presupposes periodicity, but suffices
here for our immediate purposes.

Suppose that $w \in H^s_{\ga}(\RR \times \SS^{n-1})$ for some $\ga \notin
\Gamma$ is supported in $t \ge 0$ and solves $L_{\e}w = f$ for some
compactly supported smooth function $f$ (hence $w \in H^s_{\ga}$ for
every $s$ by elliptic regularity). For example, we could
take $\phi$, which solves $L_{\e}\phi = 0$, and let $w = \chi \phi$
where $\chi$ is a cutoff function having support in $t \ge 0$ and equalling
one for $t \ge 1$. Taking transforms
we get $\hat L\hat w = \hat f$; this function on the right,
$\hat f(t,\theta,\z)$, is obviously smooth in $(t,\th)$ and
entire in $\z$. Applying the inverse $\hat G$ from the
last section gives $\hat w(t,\th,\z) = \hat G \hat f(t,\th,\z)$.
The right side of this equation is meromorphic in $\z$ with poles
at some subset of points in $\Cal P$, hence the same is true for the
left side as well. Notice that the poles of $\hat G$ in $\Im\z < -\ga$
must be cancelled by zeroes of $\hat f$ since $\hat w$ is a priori
known to be regular in this half-plane.

The function $w$ is recovered by inverting $\Cal F$,
integrating along the line $\Im\z = -\ga$.  However, by Cauchy's
theorem this
integration may be taken along any higher contour $\Im\z = -\ga'$, so
long as the interval $[\ga',\ga]$ does not contain any points of
$\Gamma$. Thus, for any such $\ga'$, $w \in H^s_{\ga'}$.
If the contour is shifted even further, so as to cross a point in
$\Cal P$, i\.e\. a pole of $\hat G$, then the resulting integral
along the line $\Im\z = -\ga''$
produces a new function $v$. Since $(L_{\e}v)\sphat =
\hat L_{\e}\hat G\hat f = \hat f$ and $\hat f$ is entire, it follows
that $v$ also solves $Lv = f$, and so
$w$ and $v$ must differ by an element of the nullspace of $L_{\e}$.
By construction $v\in H^s_{\ga''}$, and so we have decomposed
$w$ as a sum $v + \psi$, where $L_{\e}\psi = 0$. $\psi$
is given by the residues of $\hat w$ at the points $\z_j\in\Cal P$
with $-\ga < \Im(-\z_j) = -\ga_j < -\ga''$. (We use $-\z_j$ instead
of $\z_j$ for notational convenience only.) If there is just one such point,
say $-\z_j$, then it is clear that in fact
$\psi \in H^s_{\ga_j + \e}$
for every $\e > 0$. Actually, it is not hard to show, using results
from \S2, that this $\psi$ must grow or decay exactly like a polynomial
in $t$ times $e^{\ga_j t}$.

This process may be continued by moving the contour past more and
more poles. At each step, we have decomposed $w$ into a sum of
solutions of $L_{\e}\psi = 0$ and a term which decays at a rate
given by the height of the contour.
\proclaim{Proposition 4.14} If $w$ solves $L_{\e}w = f$ for
some compactly supported function $f$ on $(0,\infty) \times \SS^{n-1}$,
with $w \in H^s_{\ga}$ for some
$\ga \notin -\Gamma$, then, as $t \rightarrow \infty$,
$w(t,\th) \sim \sum \psi_j(t,\th)$, where each $\psi_j$ solves
$L_{\e}\psi_j = 0$ and $\psi_j$ decays like a polynomial in
$t$ multiplying $e^{-\ga_j t}$, where $\ga_j\in -\Gamma$ as in (4.12).
\endproclaim

Note that this result shows that the poles $\z_j$, or at least
their imaginary parts, $\ga_j$, correspond to the precise growth rates
of solutions of $L_{\e}\phi = 0$ on the cylinder.

A similar, though slightly more complicated, expansion holds for
the elements of the nullspace of the linearization $L$ in (4.2).
\proclaim{Proposition 4.15} Let $w$ solve $Lw = 0$,
at least on some neighborhood of the puncture $p_j$, and
lie in $H^s_{\ga}$ for some $\ga \notin -\Gamma$, where
$\Gamma$ is the set of imaginary parts of poles corresponding to the
inverse $\hat G$ for the model $L_{\e}$ for $L$ near this puncture. Then
$w$ has an asymptotic expansion
$$ w(t,\th) \sim \sum_{j,k = 0}^{\infty} \psi_{j,k}(t,\th)  $$
as $t \rightarrow \infty$. In this sum, the `leading terms'
$\psi_{j,0}$ are solutions of $L_{\e}\psi_{j,0}=0$ corresponding
to the poles $\z_j$ above $\Im\phi = -\ga$. These decay, as before,
like a polynomial in $t$ multiplying $e^{-\ga_j t}$. The higher terms
$\psi_{j,k}$, $k > 0$, decay like a polynomial in $t$ multiplying
$e^{-(\ga_j + k\a)t}$.
\endproclaim
\demo{Proof}
As before, we may assume that $w$ is supported in $t \ge 0$ and
$Lw = f$ is compactly supported and smooth. By (4.3),
$$
L_{\e}w = -e^{-\a t} Fw + f. \tag 4.16
$$
Now conjugate by the Fourier-Laplace transform, and apply the inverse
$\hat G$ of $\hat L_{\e}$ to get
$$
\hat w = \hat G(-e^{-\a t} Fw)\sphat + \hat G\hat f. \tag 4.17
$$
The term on the left is holomorphic, a priori, in the
half-plane $\Im\z < -\ga$. The second term on the right is
entire, while the first term on the right is holomorphic
in $\Im\z < -\ga$ and extends meromorphically to
the slightly larger half-plane $\Im\z < -\ga + \a$.
As before, $w$ can be recovered by integrating along
$\Im\z = -\ga - \e$; if this contour is moved up to $\Im\z
= -\ga + \a - \e$, and if the strip $-\ga < \Im\z < -\ga + \a$
contains no poles, then we find that $w \in H^s_{\ga -\a}$.
If there are poles in this strip,
then, exactly as in the last Proposition, $w$ decomposes
into a sum $v + \psi$ with
$v \in H^s_{\ga - \a}$ and $L_{\e}\psi = 0$. $\psi$ in turn
decomposes into a sum of terms $\psi_{j,0}$ corresponding
to the various poles in this strip.
This improved information may be then fed back
into (4.16) and (4.17). Now the right side of (4.17) extends
meromorphically to the half-plane $\Im\z < -\ga + 2\a$, and the
contour may be shifted further to get more contributions to the
expansion for $w$. Continuing this bootstrapping yields the full
expansion. Details are left to the reader.
\enddemo

\head The Deficiency Subspace \endhead

Particularly important for us in the application of
the implicit function theorem will be the pole $0 \in \Cal P$
for $\hat G$ (for any value of the Delaunay parameter $\e$).
Specifically, we will be concerned with those solutions of $L_{\e}w=0$
and $Lw=0$ which are in $H^s_{\d}$ for every $\d>0$, but not
in any $H^s_{-\d}$. Here we are still only concerned with
the local behavior of these solutions on each end $E_j$.
Their global nature will be discussed later.

The solutions to $L_{\e}\phi = 0$ in $H^s_{\d} \ominus H^s_{-\d}$
are given in (2.15) as the functions $\phi_1 = \phi_{1,\e}$ and
$\phi_2 = \phi_{2,\e}$. They depend on $t$ but not
on $\theta$. The fact that $\phi_2$ grows linearly in $t$ indicates
that $\hat G$ has a pole of order $2$ at $\phi = 0$. These
are the only temperate solutions of $L_{\e}\phi = 0$
on the whole cylinder; any other solution grows at an
exponential rate (with possible exponents given by the values
of $\Im\Cal P =\Gamma$).

If the solution metric $g$ has model Delaunay parameter $\e$
on the end $E_j$, then these functions (for that value of $\e$)
may be cut off and transplanted onto this semi-cylinder.
They do not decay, but it is easy to create a sequence of
cut-offs $\chi_i\phi$ ($\phi = \phi_1$ or $\phi_2$) with the following
properties: each $\chi_i$ is compactly supported, and
has support tending to infinity as $i$ tends to infinity.  The $L^2$ norm of
$\chi_i\phi$ equals one for all $i$, but the $L^2$ norm of $L (\chi_i\phi)$
tends to zero ($L$ is the linearization, not one of the models $L_{\e}$
here). The existence of such sequence is a standard criterion for
showing that $L$ does not have closed range on $L^2$. In any event, we
think of $\phi_1,\ \phi_2$ as constituting the bounded approximate
nullspace for $L$ (here `bounded' is loosely interpreted to encompass
the linearly growing $\phi_2$).

Define now a linear space $W$ generated by the functions
$\phi_{1,\e_j}, \phi_{2,\e_j}$, cut off in a fixed way so as to be
supported in $t \ge 0$ and transplanted on each end $E_j$; the
$\e_j$ are the Delaunay parameters for the particular model metrics
$g_{\e_j}$ on $E_j$. Since there are $k$ ends, $W$ is $2k$-dimensional.
For reasons that will become clear in the next section, we call
$W$ the deficiency subspace for $L$. Clearly $W \subset H^s_{\delta}$
for every $s$ and every $\d > 0$.

\head The Linear Decomposition Lemma \endhead

There are two important corollaries of Proposition 4.15 and
its proof which we single out in this section.
The first concerns the behavior of solutions of the inhomogeneous
equation $Lw = f$ on each end $E_j$, while the second studies the
exact value of the exponent $\a$ appearing in Theorem 3.3
and then later in (4.3), etc.

As already pointed out in Corollary 4.13, if $L$ has no
global $L^2$ nullspace, then we can find a solution $w \in
H^{s+2}_{\d}$ to the equation $Lw = f$ for every $f \in H^s_{\d}$,
whenever $\d > 0$. In particular, this holds whenever $f \in H^s_{-\d}$.
Clearly, whenever $f$ decays at some exponential rate like this,
we expect the solution $w$ to be somewhat better behaved than
a general $H^{s+2}_{\d}$ function; of course, it is immediate
that it is in this space for any $\d>0$, but we can do even better.
This is the subject of what we will call the
\proclaim{Linear Decomposition Lemma 4.18} Suppose $f \in H^s_{-\d}$ for
some $\d>0$ sufficiently small, and $w \in H^{s+2}_{\d}$ solves
$Lw = f$. Then $w \in H^{s+2}_{-\d} \oplus W$, i\.e\. $w$ may be
decomposed into a sum $v + \phi$ with $v$ decaying at the same
rate as $f$ and $\phi$ in the deficiency subspace $W$.
\endproclaim
\demo{Proof} Clearly this question may be localized to each
end $E_j$, and in this localized decomposition $\phi$ will be
a combination of $\phi_1$ and $\phi_2$. The decomposition is
achieved by exactly the same sort of shift of the contour
in the integral defining $\Cal F^{-1}$ across a pole of
$\hat G$. Here the contour is being shifted from
$\Im\z = -\d$ to $\Im\z = +\d$; the pole crossed is the one at $\z = 0$.
\enddemo

The second corollary deals with the rate at which a general
solution of (1.3) on a punctured ball, singular at $\{0\}$,
converges to the radial Delaunay metric.  We assume the
simpler statement of Theorem 3.3 that $u$ decays to
$u_{\e}$ at some exponential rate $\a$, and use the linear
theory to find the optimal rate. First transform the punctured
ball to the half-cylinder $(0,\infty) \times \SS^{n-1}$, and assume
all functions are defined here. Write the solution $u$ as
a perturbation $(1+v)u_{\e}$, so that $N_{\e}(v) = 0$, where
$N_{\e}$ is the nonlinear operator (2.11). This equation
is the same as $L_{\e}v = -Q(v)$, where $Q$ is the quadratically
vanishing function in (2.12). We already know that $v \in H^s_{-\a}$
(for all $s$), hence $Q(v) \in H^s_{-2\a}$. By the contour-shifting
arguments above, this implies that $v$ itself decays at this
faster rate, at least provided $2\a < \ga_1$, where $\ga_1$
is the first positive element in $\Gamma$. Continue this process
until this first pole at $\z_1$ has been crossed;
the conclusion is that $v$ decays exactly like $e^{-\ga_1 t}$
(possibly multiplied by a polynomial in $t$). We could
also bootstrap further and obtain
a complete asymptotic expansion for $v$, of the same general form
as in Proposition 4.15 above, although more exponents occur because of
the nonlinearity. (This is a simple form of the argument in [M2]).
We summarize this discussion as
\proclaim{Proposition 4.19} The exponent $\a$ occuring in
Theorem 3.3 governing the rate of decay of a general
solution $u$ to its model Delaunay solution $u_{\e}$ is equal to the
first nonzero element $\ga_1$ in the set $\Gamma$ corresponding
to $\hat G_{\e}$. The function $u$ admits a complete
asymptotic expansion into terms of increasingly rapid exponential
decay.
\endproclaim

\head The Bounded Nullspace \endhead

In this last subsection of \S4, we finally come to the global
behavior of that portion of the nullspace of $L$ corresponding
to the pole at $z=0$.  This space, $\BL$, which we shall call
the `bounded nullspace' of $L$, is defined by
$$
\BL = \{v \in H^s_{\d}: Lv = 0, v \notin H^s_{-\d} \}. \tag 4.20
$$
The full nullspace of $L$ on $H^s_{\d}$ will be the direct sum
of $\BL$ and the $L^2$ nullspace (which, by Proposition 4.15
is the same as the
nullspace of $L$ in $H^s_{-\d}$). In particular, when this latter space
is trivial, $\BL$ is the full nullspace of $L$ in $H^s_{\d}$.
By the linear decomposition lemma 4.18, $\BL$ is already
contained in $H^s_{-\d} \oplus W$. The purpose of this subsection is
to determine the dimension of $\BL$.

As usual, an index theorem is the principal tool for calculating
$\dim(\BL)$. Fortunately we require only a relative index theorem,
which computes the difference between two indices in terms of
asymptotic data, rather than global data.
To be more explicit, for any $\ga \notin \Gamma$ set
$$
\ind (\ga) = \dim\text{\,ker\,}\left. L\right|_{H^s_{\ga}}
- \dim\text{\,coker\,}\left. L \right|_{H^s_{\ga}}. \tag 4.21
$$
This is obviously independent of $s$. Since the adjoint of $L$ on
$H^s_{\ga}$ is $L$ on $H^{-s}_{-\ga}$, duality implies that
$$
\ind(-\ga) = -\ind(\ga)\ \text{for every\ }
\ga \notin \Gamma. \tag 4.22
$$

If $\ga_1$ and $\ga_2$ are any two allowable values
(i\.e\. neither are in $\Gamma$), the relative index
with respect to these two numbers is simply the difference
$$
\relind(\ga_1,\ga_2) = \ind(\ga_1) - \ind(\ga_2).
$$
In particular, using duality again,
$$
\aligned
\relind(\d,-\d) = & 2\,\ind(\d) \\
 = & 2\left(\dim\text{ker}\left. L \right|_{H^s_{\d}} - \dim\text{ker}\left.
L \right|_{H^s_{-\d}} \right) \\
 = & 2\dim(\BL). \endaligned  \tag 4.23 $$

As noted above, this relative index can be shown, on fairly
general principles, to be computable in terms
of asymptotic data for the operator $L$. Finding a specific
and computable formula is another matter and, to our
knowledge, there is no general result of this sort available for
asymptotically periodic operators. However, such a result is available for
for operators associated to asymptotically cylindrical metrics, and we will
use this instead. Our result is
\proclaim{Theorem 4.24} $\dim(\BL) = k$. \endproclaim
\demo{Proof}
By (4.23), it suffices to show that
$\relind(\d,-\d) = 2k$. By the usual stability properties
of the index (hence any relative index) under Fredholm deformations,
we compute this number by choosing a one-parameter family of Fredholm
operators $L^{\sigma}$, $0 \le \sigma \le 1$, with $L^{0} = L$, and
$L^1$ an operator for which this relative index is computable.
Since $L$ is just $\Lap_g + n$, we shall choose a one-parameter
family of metrics $g^{\sigma}$ and define $L^{\sigma}$ to be
$\Lap_{g^{\sigma}} + n$. The metrics $g^{\sigma}$ will agree
with $g$ except on each of the ends $E_j$, where we make the
following homotopy. First deform $g$ on each end to its model Delaunay
metric $g_{\e}$, and then deform each $g_{\e}$ through
Delaunay metrics to the cylindrical metric $g_{\bar u}$. The metric
$g^{1}$ will agree with the original $g$ on any large fixed compact
set, and will equal $g_{\bar u}$ on each end $E_j$. This metric
is now an {\it exact $b$-metric}, in the language of [Me]
(more prosaically, it has asymptotically -- and in this case, exactly --
cylindrical ends), and the corresponding operator $L^1$ is an elliptic
$b$-operator.

Of course, we still need to prove that $L^{\sigma}$ is Fredholm on
$H^s_{\d}$ and $H^s_{-\d}$ for every $0 \le \sigma \le 1$, provided
$\d>0$ is sufficiently small.  In the part of the deformation
where $g$ is homotoped to its model Delaunay metric $g_{\e}$ on each
end, this is obviously true. For the remaining part of the
deformation, we need to know that the elements of $\{\ga_j \in
\Gamma, j>1\}$, which are the weights for which $L$ will not be
Fredholm, remain bounded
away from $0$ as $\e$ varies between its initial value and $\bar u$.
This is precisely the content of Lemma 2.18.

Now, to apply Melrose's Relative Index Theorem [Me], note that we
may as well consider the operator $\tilde L^{1} = \frac {n-2}n L$,
which on each end takes the form $\del_t^2 + \del_{\th}^2 + (n-2)$.
The set $\Cal P$ for this operator (called $\text{spec}_b(\tilde L^1)$
in [Me]) is $\{\pm i\sqrt{n-2}, \pm \ga_j, j \ge 1\}$, where
$\ga_j$ are all strictly positive, tending to infinity, and
obtained in a straightforward manner from $\text{spec}(\Lap_{\th})$.
The prescription to calculate the relative index is to first consider
the `indicial family,' which in this case equals $\Lap_{\th} + (n-2)
- \zeta^2$ (it is obtained by passing to the Fourier tranform
with respect to $t$, which carries $\del_t$ to $-i\zeta$). This
is a holomorphic family of elliptic operators on $\SS^{n-1}$, and
its inverse $G(\zeta)$ has poles exactly at the points of $\Cal P$.
Since we are computing the jump in the index as the weight changes
from $-\d$ to $+\d$, we need to compute the `degree' of each pole
(as defined in [Me])
for each element of this set with imaginary part equal to zero.
Each of the poles $\pm \sqrt{n-2}$ is simple and of rank one, and
so the degree of each is also equal to one.  Finally we need to
sum over all poles {\it and} over all ends $E_j$, because the
preceding discussion is local on each end.  For each end there are
two poles of $G$ with imaginary part zero, each contributing
a degree of one to the computation, and there are $k$ ends;
The sum of all of these is $2k$, and this is the relative
index $\relind(\d,-\d)$.
\enddemo

\specialhead V. The Moduli Space: Smooth Points \endspecialhead

In this section we commence the study of the moduli space $\ML$
itself. Here we use the linear theory developed in the previous
section in a straightforward way to study
neighborhoods of the `good points' $g \in \ML$, where
the associated linearization has no $L^2$ nullspace. In the
next two sections we develop ideas to study neighborhoods
of the nonsmooth points.

The implicit function theorem is directly applicable only when
the linearized operator $L_g = \Lap_g + n$ is surjective on the
appropriate function spaces. These function spaces should be
tangent to a suitable space of metrics conformal to $g$ and
with prescribed growth conditions on each end of $\SL$.
Let $g \in \ML$, and consider all nearby metrics conformal
to $g$ of the form $(1+v)^{\frac 4{n-2}}g$. The most natural
class of such metrics is
$$
\text{Met}^{s+2}_{-\d} = \{(1+v)^{\frac 4{n-2}}g, \ v
\in H^s_{-\d}, \ v > -1\}. \tag 5.1
$$
Clearly $T_g\text{Met}_{-\d}^{s+2} = H^{s+2}_{-\d}$. The operator
$N_g(v)$, as defined in (4.1), can be considered as a map
from this space of metrics to $H^s_{-\d}$, and as such
is obviously $\cin$. However, it is never surjective, since $-\d < 0$.
The obvious alternative, to consider metrics growing at the rate $+\d$
is unsuitable because of the nonlinearity.

The case where this difficulty is easiest to remedy is
when $L_g$ has no $L^2$ nullspace:
$$ \text{ker\,}(L_g) \cap L^2(\SL,dV_g) = \{0\}. \tag 5.2 $$
This will be our standing hypothesis in the rest of this section.

We can rephrase this condition by recalling, by virtue of Proposition
4.15 on asymptotics of solutions of $L\phi = 0$, that if $\phi$ were in
the $L^2$ nullspace of $L$, then it would decay at the exponential rate
$e^{-\ga_1 t}$ on each end, where $\ga_1 = \ga_1(\e_j)$ is the first
nonzero element of $\Cal P_{\e_j}$ on the end $E_j$, corresponding to
the model Delaunay metric $g_{\e_j}$ there. In particular,
under the hypothesis (5.2), for $\d \in (0,\ga_1)$, the operator $L$
has no nullspace on $H^s_{-\d}$ for any $s \in \RR$. Now apply
Corollary 4.13 and the linear decomposition lemma 4.18
to conclude that (5.2) is equivalent to the statement that
$$
\gathered L: H^{s+2}_{-\d}(\SL, dV_g) \oplus W \longrightarrow
H^s_{-d}(\SL,dV_g)\\
\text{is surjective for any\ } 0 < \d < \min_j\{\ga_1(\e_j)\}
\ \text{and\ } s \in \RR. \endgathered \tag 5.3
$$
$W$ is the $2k$-dimensional deficiency
subspace introduced in the last section. Because we work with a nonlinear
equation, we assume that $s > n/2$.

\proclaim{Theorem 5.4} Suppose $g \in \ML$ is a metric such that
the hypothesis (5.2) is satisfied. Then there exists a $2k$-dimensional
open manifold $\Bbb W$, the elements of which are functions $v$ on
$\SL$ to be described below, and with a distinguished element
$v_0 \in \Bbb W$, such that the nonlinear map
$$
N_g: H^{s+2}_{-\d}(\SL; dV_g) \times \Bbb W \longrightarrow
H^s_{-\d}(\SL; dV_g)
$$
defined by $N_g(v,\psi) = N_g(v+\psi)$ where $N_g$ is the
operator given by (4.1),
restricted to a neighborhood of $0$ in the first factor, is
real analytic. Moreover, the tangent space $T_{v_0}\Bbb W$ is
identified with the deficiency subspace $W$, and using this identification,
the linearization $L_g$ of the map $N_g$ is surjective, as in
(5.3).
\endproclaim
\proclaim{Corollary 5.5} Suppose $g\in\ML$ satisfies the hypothesis (5.2).
Then there is an open set $\Cal U \subset \ML$ containing
$g$, such that $\Cal U$ is a $k$-dimensional real analytic manifold.
\endproclaim

The proof of the Corollary follows directly from the surjectivity
statement (5.3), by a straightforward application of the real
analytic implicit function theorem.
The dimension $k$ is of course
the dimension of the nullspace of $L$ as a map (5.3). This nullspace
is the bounded nullspace $\BL$ defined in the last section, and
by Theorem 4.24 its dimension is $k$. So to conclude the proof of
the theorem and the corollary it suffices to construct $\Bbb W$, and
to show that $N_g$ is a real analytic mapping.

The only difficulty in the construction of $\Bbb W$ is that
the elements of $W$ are not, generally, bounded either above or below,
so we could not
use $(1+v)^{4/(n-2)}$ as a conformal factor for most $v \in W$.
However, recall from (2.15) that each element $\phi_{j,\e}$, $j = 1,2$,
of the bounded nullspace for the model problem on the cylinder is the
tangent to a curve $\eta \rightarrow \Phi_{j}(t,\e,\eta)$ of actual
conformal factors. The actual curve $v_{\eta} \in \Bbb W$ we would want to
use for each such element of $W$ on an end has the form
$$
v_{\eta}(t) = \cases \frac{u_{\e}(t + \eta)}{u_{\e}(t)} - 1, &j = 1 \\
\frac{u_{\e + \eta}(t)}{u_{\e}(t)} - 1, &j = 2. \endcases
$$
The manifold $\Bbb W$ is constructed by gluing together these
local definitions from each end $E_j$ in an essentially arbitrary,
but smooth, manner. By construction, $T_0\Bbb W = W$, as desired.

To check the real analyticity of $N_g$, the main step is to write
$N_g(v + \psi)$ in such a way that it clearly lands
in $H^s_{-\d}$. Hence, on each end write $N_g = N_{\e} + e^{-\e t}Q$,
and then use a common formula for the remainder in Taylor's theorem
to obtain
$$ N_g(v+ \psi) = N_{\e}(\psi) + e^{-\e t}Q(v+\psi) + \left[
\int_0^1 N_g'(\psi + sv)\,ds\right]v. $$
Since $N_{\e}(\psi) = 0$, every term on the right here is in $H^s_{\d}$
so long as $\d < \e$. Furthermore, every term is real analytic
in $(v,\psi)$. Real analyticity of $N_g$ on the interior, away
from the ends, is even easier.

The rest of the proof is now standard.

\specialhead VI. The Real Analytic Structure of $\ML$ \endspecialhead

In this section we prove that $\ML$ is always an analytic set
by representing it as the slice of an infinite dimensional
real analytic manifold $\MML$ with the conformal class $g_0$.
We go on to prove a generic slice theorem, that even if $\ML$ itself
is not smooth, generic nearby slices of $\MML$ by other conformal
classes will be.

\head The Urmoduli Space $\MML$ \endhead

As described above, we wish to regard $\ML$ as the slice by
the standard conformal class $[g_0]$ of the infinite dimensional
set $\MML$ consisting of all metrics, not necessarily in
the standard conformal class on $\SS^n$, with constant scalar
curvature $n(n-1)$. For obvious technical reasons, we
consider only those metrics which satisfy appropriate
growth and asymptotics conditions near the punctures $p_j$
(this will be be elaborated on below.)
Consideration of this `big,' or `ur-,' moduli space (when the
underlying manifold is compact) was
first undertaken by Fischer and Marsden in the early
'70's ([FM1] and [FM2]), motivated in part by concerns in general
relativity.  They proved that it is a smooth Banach manifold (in any
one of of a number of standard Banach completions), provided a certain
overdetermined linear equation has only the trivial
solution.

Their precise set-up was to consider $\Bbb M_{\rho}$, the set
of all metrics (of some fixed finite regularity) with scalar
curvature function $\rho$, a fixed function on the underlying
manifold $X$. To show that this set is a Banach manifold, it
suffices to show that the linearization of the scalar curvature
map (which assigns to any metric $g$ its scalar curvature
function $R(g)$) is surjective at any $g \in \Bbb M_{\rho}$.
Since this map $R$ carries metrics to functions, its linearization,
which we will call $\a_g$ here, carries the tangent space at $g$
of the space of all metrics, i\.e\. the symmetric $2$-tensors,
to scalar functions. A. Lichnerowicz [Li] had earlier computed that,
for $h$ a symmetric $2$-tensor,
$$ \a_g(h) = \d\d h - \Lap_g(\text{tr\,}h) - \langle r_g, h \rangle.
\tag 6.1 $$
Here $\d$ is the divergence operator on tensors, and $r_g$ is
the Ricci tensor for $g$. It is easy to check that $\a_g$ reduces
to (a multiple of) the linearization operator $L_g$ we have
been studying when $h$ is a multiple of the metric $g$ (and hence is
tangent to the conformal class $[g]$), when $\rho = n(n-1)$.
It is straightforward
to check that the symbol of $\a_g$ is surjective, so it follows
that the symbol of its adjoint $\a_g^*$ is injective, and then
that the symbol of $\a_g \circ \a_g^*$ is an isomorphism. Therefore,
this last, fourth order operator is elliptic; at least when $X$ is
compact, it is Fredholm, and in particular has closed range of finite
codimension. Since $\text{range\,}(\a_g) \supset \text{range\,}(\a_g
\circ\a_g^*)$, the range of $\a_g$ itself is closed and of finite
codimension. (It does, however, always have an infinite dimensional kernel.)
By virtue of all this, when $X$ is compact, $\a_g$ is surjective if and only
if $\a_g^*$ is injective, and this ($\a_g^*(f) = 0$) is precisely the
overdetermined equation referred to earlier.

{}From (6.1) it follows readily that
$$
\a_g^*(f) = \text{Hess\,}_g(f) - g \Lap_{g}(f) - f\,r_g. \tag 6.2
$$
If $\a_g^*(f) = 0$, then taking the trace we
also get that $\Lap_g(f) = -(\rho/(n-1))f$; substituting this back
into (6.2) we obtain finally the condition
$$
\a_g^*(f) = 0 \Longrightarrow \text{Hess\,}_g(f) = -\left(
r_{g} - \frac{\rho}{n-1}g \right) f. \tag 6.3
$$
Fischer and Marsden were able to show that (6.3) (or equivalently, (6.2))
has only the
trivial solution $f=0$ except possibly in the special cases where
$\rho$ is a nonnegative constant and
$\rho/(n-1)$ is in the spectrum of $-\Lap_g$.
Thus, except for these cases,
$\Bbb M_{\rho}$ is a smooth Banach manifold.

In the case of interest here, the underlying manifold $\SL$ is not
compact, and $\rho/(n-1) = n$ is {\it always} in $\text{spec\,}(\Lap_g)$.
We address the former of these concerns first, since
the Fredholm theory of \S4 may be adapted easily to this problem.
As before $H^s_{\ga}(\SL)$ will denote the Sobolev space of
scalar functions decaying like $e^{\ga t}$ on each end,
with respect to fixed cylindrical coordinates $(t,\th)$ there and with
respect to a fixed metric $g \in \ML$ (any element of the moduli space
can be used as a fixed background metric to define these Sobolev
spaces).  Define also $H^s_{\ga}(\SL, \text{S}^2T^*(\SL))$ to
be the space of symmetric $2$-tensors $h = e^{\ga t}k$, with
$k$ in the unweighted Sobolev space $H^s(\SL, \text{S}^2 T^*(\SL);
dV_g)$.
\proclaim{Proposition 6.4} For any $g \in \ML$, the operator
(6.1) is bounded as a linear map
$$
\a_g: H^{s+2}_{\ga}(\SL, \text{Sym}^2(T^*(\SL))) \longrightarrow
H^s_{\ga}(\SL).
$$
It has an infinite dimensional nullspace, and the closure
of its range has finite codimension for all $\ga \in \RR$.
There exists a discrete set $\Gamma'$ such that for all
$\ga \notin \Gamma'$, this map also has closed range.
\endproclaim
\demo{Proof} The boundedness assertion is immediate.
The proof of the closed range part is very close
to the analogous one given in \S4 for the scalar operator
$L_g$, so we just provide a sketch. As in [FM2] and above, the
surjectivity of the symbol of $\a_g$ guarantees that we
merely have to show that $A_g = \a_g \circ \a_g^*$, which
is elliptic, has closed range for all $\ga \notin \Gamma'$.
The operator $A_g$ is asymptotically periodic on each end
$E_j$ of $\SL$, and the Fourier-Laplace transform may again
be employed to construct local parametrices on each of these
ends. The fact needed to make this work is that the holomorphic
family of elliptic operators $\tilde A_{\e}(\phi)$, on the compact
manifold $\SS^1 \times \SS^{n-1}$, and associated to the
Delaunay model operator $A_{g_\e}$, is invertible at some
$\phi \in \Bbb C$. In this case
the analytic Fredholm theory, and the rest of the construction
in \S4, proceeds exactly as before. Again we prove this for
$\phi \in \RR$, since for such $\phi$, $\tilde A_{\e}(\phi)$
is self-adjoint. Invertibility is equivalent to the injectivity
of the induced operator $\tilde \a_{\e}^*(\phi)$, and this in turn
is implied (by taking traces) by the injectivity of
the model operator $\tilde L_{\e}(\phi)$. The invertibility of
this last operator for all $\phi \in (0,2\pi)$ (we are again
assuming the period of $\a_{\e}$, et al\., to be one for simplicity)
was already proved in \S4, and so we are done.
\enddemo

To proceed with the program of [FM2], we also have to show that
$\a_g^*$ is injective on $\SL$ for any $g \in \ML$. If this
were true, then the urmoduli space $\MML$ would be a real analytic
manifold in a neighborhood of every $g \in \ML$, i\.e\. there
would exist a (local) thickening of $\ML$ into the real analytic
moduli space $\MML$ such that $\ML = \MML \cap [g_0]$, as desired.
Unfortunately this is not quite so simple, since, as already
noted, $\rho/(n-1)$ is always in the spectrum of $\Lap_g$; we
are thus precisely in the case not treated by [FM2].

The equation (6.3) is rather similar to the well-known equation of
Obata, which also appears in the study of
the scalar curvature, and actually reduces to Obata's equation
when $X$ is an Einstein manifold. In analogy with the situation there --
in particular, Obata's characterization [0] of $(\SS^n,g_0)$
as the only compact
manifold admitting nontrivial solutions of his equation --  Fischer
and Marsden conjectured that (6.3) admits a nontrivial solution, at
least when $X$ is compact, if and only if $X$ is the standard sphere.
This is now known to be false: indeed, let $X$ be the product
of the circle of length $L$, $\SS^1(L)$, with any Einstein
manifold $E^{n-1}$ of positive scalar curvature, e\.g\.
$\SS^{n-1}$. Then, if the scalar curvature of $E$ is normalized to
be $(n-1)(n-2)$, (6.3) implies that $f$ depends only
on the variable $t$ along the circle, and satisfies the ODE
$f'' + (n-2)f = 0$. Thus, (6.3) admits a nontrivial solution
whenever $L$ is an integer multiple of $\sqrt{n-2}$.
There are also noncompact counterexamples to their conjecture:
the simplest are just the universal covers $\RR \times E^{n-1}$
of the compact examples above.

The characterization even of all compact manifolds admitting
solutions to (6.3) is still unknown, although it is almost
certain that the list of possibilities is rather small.
Fortunately, though, the characterization of all complete
conformally flat manifolds admitting nontrivial solutions
to (6.3) was obtained about ten years ago,
independently by J. Lafontaine [L] and O\. Kobayashi [Ko].
The simplification in this case is
that the Ricci tensor of any conformally flat manifold is
always harmonic. Kobayashi and Lafontaine proved that beyond
the sphere and the product type examples discussed above, there
are also a collection of examples which are warped products
of $\SS^1$ with $E^{n-1}$. We shall not write these examples
explicitly here, but simply note that, provided $k > 2$, the
solution metrics $(\SL,g)$, $g \in \ML$, are never warped
products (for topological reasons alone!). We obtain then
\proclaim{Proposition 6.5} There are no nontrivial solutions to
the equation $\a_g^*(f)$ for any $g \in \ML$, provided $k>2$.
\endproclaim
Even when $k=2$, only the cylinder $\RR \times \SS^{n-1}$ amongst
the Delaunay metrics appears on the list of counterexamples.
To check this, use the system of equations given in Lemma 1.1 of [Ko].
These involve both the warping function for
the metric and the function $f$. Elementary
manipulations show that the noncylindrical Delaunay metrics
never satisfy this system.
\proclaim{Corollary 6.6} $\MML$ is a real analytic
Banach manifold in a neighborhood of $\ML$.
\endproclaim
\demo{Proof}
It remains only to set up the precise function spaces on
which the implicit function theorem (see [F] page 239)
will be applied, and make a few additional comments.
We use, of course, the weighted
Sobolev spaces $H^s_{\ga}$ of symmetric $2$-tensors and functions,
as before.

The crux of the argument is the fact that when $g \in \ML$,
$\a_g$ is surjective both on $H^{s}_{\d}$ {\it and} on $H^{s}
_{-\d}$, and indeed on $H^{s}_\ga$ for any $\ga \notin \Gamma'$!
This seems somewhat counterintuitive, but follows from
the highly overdetermined nature of $\a_g^*$. In fact, we already
know that $\a_g$ has closed range on $H^{s}_\ga$, $\ga \notin \Gamma'$
by Proposition 6.4. Its cokernel is identified with the kernel of
the adjoint, $\a_g^*$ which is a map from $H^{-s}_{-\ga}$ to
$H^{-s-2}_{-\ga}$. By injectivity of the symbol of $\a_g^*$
we may use elliptic regularity to replace $-s$ here by
any positive number, e\.g\. $+s$. Finally we can invoke
the results of Lafontaine and Kobayashi, via Proposition 6.5,
to conclude that $\a_g^*$ has no nullspace, regardless of the
weight $-\ga$.

Now consider the scalar curvature map
$$ \Cal N_g: H^{s+2}_{-\d}(\text{Sym}^2) \times \Bbb W
\longrightarrow H^s_{-\d}(\text{Sym}^2), \tag 6.7 $$
defined in the obvious way. We are including the factor $\Bbb W$
here so that points of the urmoduli space are allowed to have
varying neck sizes, but all other perturbations are required to
decrease exponentially on the ends. To apply the analytic implicit
function theorem we need to know, first, that the linearization
$$ \a_g: H^{s+2}_{-\d}(\text{Sym}^2) \oplus W \longrightarrow
H^s_{-\d} \tag 6.8$$
is surjective, which we have just established, and secondly
that (6.7) is a real analytic mapping of Hilbert spaces.
This latter statement is also straightforward, so this completes the proof.
\enddemo

One unfortunate shortfall of this theorem is that, although $\MML$
provides a thickening of $\ML$, it is not a uniform thickening, i\.e\.
one of fixed `width' around $\ML$. However, over any compact
subset of $\ML$ we can ensure the existence of $\MML$ out to some
fixed distance by the obvious covering argument. To find this fixed width
thickening uniformly on $\ML$ we would need to understand more about
the compactification $\CML$. If the analogues of the result of Kobayashi
and Lafontaine as well as the asymptotics Theorem 3.3 and compactness
result, Proposition 3.12, were known for non-conformally flat metrics,
this would be unnecessary, and many of the results here could be given
a more satisfactory form.

\head Analyticity of $\ML$ and the Generic Slice Theorem \endhead

In addition to its intrinsic nature and relationship to the
concerns of this paper, the urmoduli space $\MML$ is required
for two purposes: to show that $\ML$ itself is analytic, and to prove
the generic slice theorem. The proofs of these results are very
closely related, so we develop the preliminaries for them
simultaneously.

Inside the ambient space $\Cal V = H^{s+2}_{-\d}(\text{Sym}^2) \times
\Bbb W$ are two analytic (Hilbert) submanifolds, namely
$\MML$ and the conformal class $[g_0]$. (By $[g_0]$ we always
mean the set of all metrics in $\Cal V$ conformal to the round
metric $g_0$.)
Likewise we have all other conformal classes in $\Cal V$; a generic
such class will be denoted $[g']$. Notice that all the analytic
machinery, in particular the Fredholm theory and relative
index computations in \S4, still holds
for any metric (corresponding to an element) in $\Cal V$.

Choose any element $g \in \MML$ and let $[g']$ be the conformal
class of $g$. We shall assume that $[g']$ is very near
$[g_0]$. The immediate aim is to prove that $\MML$ and $[g']$
intersect almost transversely at $g$; in general these submanifolds
may not be transverse at $g$, however they will always form
a Fredholm pair there. This means that their tangent spaces,
$E \equiv T_g[g']$ and $F \equiv T_g\MML$, which are closed
linear subspaces in $V \equiv T_g\Cal V = H^{s+2}_{-\d}
\oplus W$ form a Fredholm pair,
i\.e\. that $E \cap F$ is finite dimensional, and that $E+F$ is closed
in $V$ and has finite codimension there. To any Fredholm pair
one can associate an index, which is the dimension of $E \cap F$
minus the dimension of $V/(E+F)$. This number is stable under
perturbations of the pair. These matters are explained more
thoroughly in [K].

\proclaim{Lemma 6.9} $\MML$
and $[g']$ are a Fredholm pair at any point $g$ in their
intersection.
\endproclaim
\demo{Proof} The orthogonal complements of $E$ and $F$ are
given by
$$ \aligned
E^{\perp} &= \{h: \tr_g(h) = 0\}, \\
F^{\perp} &= \{h: h = \a_g^*(f), \ f \in H^{s+4}_{-\d}\}.
\endaligned \tag 6.10
$$
It will suffice to show that $E\cap F$ and $E^{\perp}\cap F^{\perp}$
are finite dimensional, and that the orthogonal projection
$\pi: F \rightarrow E^{\perp}$ has closed range.

First suppose that $h \in E\cap F$. Then $h = f\cdot g$
with $f \in H^{s+2}_{-\d} \oplus W$ and $\a_g(h) = 0$.
Computing $\a_g(f\cdot g)$ we find that $L_g(f) = 0$. Since $L$ is
a Fredholm operator, there can be at most a finite dimensional
family of such $h$ in the intersection; the precise dimension
is the same as the nullspace of $L$ on $H^{s+2}_{\d}$.
On the other hand, if $h \in E^{\perp} \cap F^{\perp}$,
then $h = \a_g^*(f)$ with $f \in H^{s+4}_{-\d}$ and $\tr_g(h) = 0$.
But $\tr_g(\a_g^*(f))$ is, up to a factor, just $L_g(f)$, so
the dimension of this intersection is the same as the dimension
of the nullspace of $L$ on $H^{s+2}_{-\d}$, and is thus
finite.

To conclude the proof, we need to know that the projection
$\pi$ has closed range (and hence is Fredholm, by the work above).
The map $\pi$ is defined, for any $h \in V$, by decomposing
$h = \lambda\cdot g + k$, where $\tr(k) = 0$, and setting $\pi(h) = k$.
Now suppose $h_j$ is a sequence of elements in $F$, $\|h_j\| = 1$
for all $j$, but $k_j = \pi(h_j) \rightarrow 0$ in norm.
Since $\a_g(h_j) = 0$ we have $(n-1)L_g(\lambda_j) =
\d \d k_j - \langle r_g, k_j \rangle$. The right hand side of
this equation certainly goes to zero in norm, so by our
Fredholm theory for $L_g$, the functions $\lambda_j$
decompose as $\lambda_j = \mu_j + \nu_j$ with $\nu_j$
going to zero in norm and $\mu_j$ in the nullspace.
Since $\mu_j \cdot g$ is in the nullspace of $\a_g$ also,
we may subtract it off from $h_j$, and get that $\|h_j\| \rightarrow 0$,
contrary to hypothesis. Hence $\pi$ restricted to $F$ has
closed range, and this completes the proof.
\enddemo

Before stating the next result we need to introduce some notation.
First, if $[g']$ is a conformal class on $\SS^n$ (represented by elements
in $\Cal V$, and assumed near to $[g_0]$), then $\ML([g'])$ will denote
the moduli
space of complete metrics on $\SL$ with constant scalar curvature $n(n-1)$
in the conformal class $[g']$. For $\e>0$ any sufficiently small number
we also let $\Cal M_{\Lambda, \e}([g'])$ denote the subset of $\ML([g'])$
consisting of solution metrics $g$ with the
Killing norms of all Poho\v{z}aev invariants bounded below by $\e$
(see \S7 for a discussion of these norms).
This is simply the subset of $\ML([g'])$ consisting of solutions
with necksizes uniformly bounded away from $0$ by $\e$.
By Proposition 3.12, each $\Cal M_{\Lambda,\e}([g_0])$ is compact;
then, provided $[g']$ is sufficiently close to $[g_0]$,
$\Cal M_{\Lambda,\e}([g'])$ will also be compact, since it is contained
in a compact neighborhood of $\Cal M_{\Lambda,\e}$ by virtue of the
implicit function theorem and the considerations above.

\proclaim{Generic Slice Theorem 6.11} For any fixed $\e>0$,
the truncated moduli space $\Cal M_{\Lambda,\e}([g'])$
is a $k$-dimensional real analytic manifold for all conformal classes $[g']$
in a set of second category and sufficiently near to $[g_0]$.
\endproclaim
\demo{Proof} The reason for introducing the truncated moduli spaces is,
of course, because we only know the existence of $\MML$ to
a fixed distance away from $\ML$ only over each $\Cal M_{\Lambda,\e}$.
We shall not comment further on this modification, but simply
remark that it could be removed provided somewhat more were known
about $\MML$, as discussed at the end of the previous subsection.

This result follows from the Sard-Smale
theorem, once we have checked the hypotheses. The basic point,
of course, is that $\ML([g'])$ is the same as $\MML \cap [g']$.
We parametrize the set of conformal classes (not modulo
diffeomorphisms!) near to $[g_0]$
by the linear Hilbert space $E^{\perp} =
\{k \in H^{s+2}_{-\d}: \tr_g(k) = 0\}$ for
some fixed $g \in \ML$. Now consider the projection map
$$
\Pi: \MML \longrightarrow E^{\perp}. \tag 6.12
$$
$\ML$ is given by the preimage $\Pi^{-1}(0)$. The preimage
$\Pi^{-1}(k)$ in $\MML$ is precisely the set $\ML([g'])$
where $[g']=[g_0+ k]$ is the conformal class corresponding to $k$.
This preimage is a smooth, in fact real analytic, manifold
provided $k$ is a regular value of $\Pi$. By Sard-Smale,
once we know that $\Pi$ is a Fredholm map (of index $k$),
then the set of regular values in $E^{\perp}$ is of second category.
Observing that the tangents to the orbits of the diffeomorphism
group are certainly not of second category, we may conclude
that the moduli spaces over generic conformal classes near to $[g_0]$
are smooth and $k$ dimensional.

However, the assertion that $\Pi$ is Fredholm of index $k$ is
contained within the preceding Lemma. Indeed,
if $g \in \MML$ and $F = T_g(\MML)$, then we proved there
that the projection $\pi = \Pi_*:F \rightarrow E^{\perp}$ is Fredholm,
and its index is the same as the relative index of $L_g$ across
the weight $0$, i\.e\. equal to $k$.

To finish the proof we need to eliminate the possibility that
$\Pi(\MML)$ is contained within a (finite codimensional) submanifold of
$E^{\perp}$. Although this would follow from knowing that $\pi$
is surjective, this is a difficult issue. Instead we appeal
to the existence theory for the nonlinear equation, generalizing
Schoen's basic construction. In fact, the construction of solutions
given in [S2] may be carried out not only for
the standard round metric $g_0$, but also for generic metrics
$g'$ which are small perturbations of $g_0$ compactly
supported away from $\Lambda$. The only modification of Schoen's
method needed to accomplish this is given in [P1], Proposition 1.2.
The linearization of this set of compactly supported perturbations is
clearly dense in $E^{\perp}$,
hence there is a dense set of $k \in E^{\perp}$ near the origin
for which the preimage $\Pi^{-1}(k)$ is nonempty. This
is sufficient to conclude that $\Pi(\MML)$ contains a full
neighborhood of $0$, hence the condition `regular value'
for $\Pi$ is not the empty one.
This completes the proof.
\enddemo

The reader should note that an alternate possibility for
proving smoothness of generic slices would be to study
the linear operator $L_g$ for any $g \in \MML$ near to
but not in the standard conformal class and show that
this operator generically has no $L^2$ kernel. The smoothness would
then follow by the results of \S5, combined with a compactness
argument. However, proving that $L_g$ has no decaying
eigenfunctions would require setting up a somewhat elaborate
perturbation theory, since this is equivalent to trying to
perturb point spectrum which is sitting on the end of a band
of continuous spectrum. Thankfully we have been able to
avoid this approach here (although we have had to appeal to the rather
more difficult existence theory for the nonlinear equation
instead). It should be noted that the surjectivity of $\Pi_*$
(which we have established) is equivalent to this eigenvalue
perturbation result.

The second main result, that $\ML$ is (locally) a real
analytic variety, also follows from Lemma 6.9.
\proclaim{Theorem 6.13} For any $g \in \ML$ there exists a
ball $B$ in $\Cal V = H^{s+2}_{-\d}(\SL; \text{Sym}^2) \times \Bbb W$,
a finite dimensional space $M$, a real analytic variety
$\Cal A$ in $M$  and a real analytic diffeomorphism
$$
\Phi: B \longrightarrow \Cal V
$$
such that $\Phi(\ML \cap B) = \Cal A \cap B'$,
where $B'$ is a small ball containing $\Phi(g)$.
\endproclaim
The proof reduces to the following abstract result, which
is presumably well-known:
\proclaim{Lemma 6.14} Let $V$ be a Hilbert space, and
$\Cal E$ and $\Cal F$ two real analytic submanifolds, the tangent
spaces of which at any point of the intersection $p \in \Gamma
\equiv \Cal E \cap \Cal F$ form a Fredholm pair. Then
for each point $p \in \Gamma$ there exists a neighborhood
$B$ of $p$, a finite dimensional subspace
$M \subset V$, a real analytic variety $\Cal A \subset M$,
and a real analytic diffeomorphism
$$
\Phi: B \longrightarrow V
$$
such that $\Phi(\Gamma \cap B) = \Cal A \cap B' \subset M \cap B'$,
where $B'$ is a neighborhood of $\Phi(p)$ in $M$.
\endproclaim
\demo{Proof} By initially composing with a real analytic diffeomorphism
we may assume that a neighborhood of $p$ in $\Cal E$ lies in a linear
subspace and that $p = 0$. (In our situation, $\Cal E$ is the conformal
class and this may be
accomplished by a translation.) Let $\bar \Psi$ denote a `defining
function' for $\Cal F$ in $V$. By this we mean that $\bar \Psi$ is a
real analytic map
from $V$ into another Hilbert space $U$ such that $\Cal F = \bar \Psi^{-1}(0)$.
By the stability of Fredholm pairs, any level set $\bar \Psi^{-1}(w)$
forms a Fredholm pair with $\Cal E$, for all sufficiently small
$w \in U$, and for all points of intersection near the origin in $\Cal E$.
It is evident that this statement is equivalent to the assertion
that $\Psi \equiv \left. \bar \Psi \right|_{\Cal E}: \Cal E \rightarrow U$
is a Fredholm map. The intersection $\Cal E \cap \Cal F$ is the same
as $\Psi^{-1}(0)$. Thus the lemma may be rephrased as saying
that if $\Psi$ is a real analytic Fredholm map between two
Hilbert spaces $\Cal E$ and $U$, then the level set $\Psi^{-1}(0)$
is locally equivalent, by an analytic diffeomorphism, to a finite
dimensional analytic variety.

A standard implicit function theorem argument proves this assertion.
Let $A = \left. \Psi_* \right|_0$, and set $M = \ker(A)$, $L =
W \ominus \text{range}(A)$. Then both $L$ and $M$ are finite
dimensional. Define a new map $\tilde \Psi: \Cal E \oplus L
\rightarrow U$ by $\tilde \Psi(q,f) = \Psi(q) + f$. Then
$$
\left. \tilde \Psi_* \right|_{(0,0)}(h,\phi) = A(h) + \phi,
$$
so this differential is obviously surjective. Note that
$$
\ker(\left. \tilde \Psi_* \right|_{(0,0)}) =
\{(h,\phi): A(h) + \phi = 0\}.
$$
Since $A(h)$ and $\phi$ lie in
orthogonal spaces, this nullspace is just $M = \ker(A)$.
The analytic implicit function theorem gives the existence
of two analytic maps
$$
k: M \rightarrow M^{\perp} \subset \Cal E, \qquad \ell: M \rightarrow L
$$
and a ball $\tilde B$ around $(0,0) \subset \Cal E \oplus L$ such
that all zeroes of $\tilde \Psi$ in this ball lie in
the graph of the map $(k,\ell): M \rightarrow M^{\perp} \oplus L$:
$$
\{(q,f): \Psi(q) + f = 0\} \cap \tilde B = \{(m + k(m),\ell(m)):
m \in M \cap \tilde B\}.
$$
Thus, for $m$ in this ball, $\Psi(m + k(m)) + \ell(m)$
vanishes identically. The zeroes of $\Psi$ in $\tilde B \cap \Cal E$
are then just $\{m+k(m): m \in \tilde B \cap M: \ell(m) = 0\}$,
for these points are all zeroes of $\tilde \Psi$, and are the
only such zeroes at which $\Psi$ also vanishes. If $\Pi$ is
an analytic diffeomorphism of $\Cal E$ carrying the graph of
$k$ into $M$, then $B \cap \Psi^{-1}(0)$ equals
$\Pi^{-1}(M \cap \tilde B) \cap \{m \in M \cap \tilde B: \ell(m) = 0\}$
as desired. This completes the proof of Lemma 6.14.
\enddemo

\specialhead VII. Concluding Remarks and Informed Speculation \endspecialhead

\head $L^2$ Nullspace and Singularities of $\ML$ \endhead

We have not yet discussed whether it is possible to give
conditions ensuring that a given $g \in \ML$ is a `smooth point,'
as defined in \S5. It is possible for $\ML$ to be both smooth
and of the correct dimension near $g$ even if $L_g$ has
$L^{2}-$nullspace; however, absence of this nullspace is our
only criterion for guaranteeing this smoothness.
We expect that it should be difficult to establish such a criterion
in general.  It would also be quite interesting to understand when $\ML$
is not smooth in the neighborhood of some element $g$. In particular,
constructing solutions near which $\ML$ is singular seems like
another very challenging and important problem.

On the other hand, the two known constructions for producing points in
$\ML$, those of Schoen [S2] and [MPU], yield solutions with explicit
geometries. The dipole solutions of [MPU], described
below, are manifestly
smooth points. These solutions exist only for certain configurations
$\Lambda$, whereas the ones in [S2] exist for any $\Lambda$.
Unfortunately, the existence of $L^2-$nullspace for these
latter solutions is less evident, although we expect this to be true;
outlined below is a strategy to prove this.
Note that once the existence of one smooth
point in any component of $\ML$ is established, the real analyticity
of the moduli space then shows that almost every element in that component
is a smooth point.

We now give a brief description of Schoen's solutions and their construction.
These solutions are uniformly small $\Cal C^{0}$ perturbations of explicit
approximate solutions. Each of these
approximate solutions is constructed from an `admissible conformal
structure' $(\Gamma,\s)$;  $\Gamma$ is an infinite tree
with a labeling, $\sigma$, of strong
dilations $G_{ij}$ for each directed edge $e_{ij}$.
We assume for simplicity that $\Gamma$ has one vertex of order $k$ and
all other vertices of order 2.  The admissibility
of the labeling $\s$ refers to certain compatibility conditions that the
dilations must satisfy.  The strength, $\lambda_{ij}$, of these dilations can
be related to a measurement of the `neck sizes' $\e_{ij}$ of the approximate
solution $g_{\sigma}$ constructed from the data $(\Gamma,\sigma)$. $(\SL,
g_{\sigma})$ consists of almost spherical regions, corresponding to
the vertices $V$ of $\Gamma$, joined by
small necks, corresponding to the edges, the sizes of which are dictated
by the $\e_{ij}$.
There are an infinite number of parameters in the construction of these
approximate solutions since the admissible
conformal structures $\s$ can be varied, even with $\Gamma$
fixed.  In particular, one
can begin with an initial conformal structure $\sigma_{1}$,
so that the corresponding approximate solution $g_{1}$ consists of one
central spherical region $\Omega_0$ (the vertex of order $k$) and $k$
periodic,
spherically symmetric ends.  But in order to find an exact solution,
it is necessary to deform the conformal structure, as dictated by the
Poho\v zaev balancing condition, and break this symmetry.

The conformal structure $(\Gamma,\sigma)$ decomposes $\SL$ into a union
of almost spherical regions, $\Omega_{i}$, $i \in V$. Each $\Omega_i
\subset \Sn$ is
the pullback by a conformal diffeomorphism $F_i$ of $\Sn$ of a
large region of $\Sn$.  The metric $g_{\s}$ on $\Omega_i$
is constructed so that $F_i: (\Omega_i,g_{\s})
\rightarrow (\Sn,g_0)$ is an isometry off of a small neighborhood
of $\del \Omega_i$. This decomposition gives rise to a basic analytic
property of the approximate solution metric $g_{\sigma}$. This is
the existence
of a basis of functions, $\eta^{\s}_{ij}$, $i \in V$ and $j = 1, \cdots, n+1$,
for an infinite dimensional space $K$ (the `small eigenspace')
corresponding to all the spectrum
of $L_{g_{\sigma}}=\Lap_{g_{\sigma}} + n$ in a small interval around $0$.
The functions $\eta^{\s}_{ij}$ have explicit geometric
descriptions. Each $\eta^{\s}_{ij}$ has
support concentrated on $\Omega_{i}$ and corresponds to
the linear function $\eta_{j}$ on $F_{i}(\Omega_{i})\subset \Sn$.
Linear functions are of course eigenfunctions for $\Lap_{g_0}$
with eigenvalue $n$. On the orthogonal complement
$K^{\perp}$, $L_{g_{\sigma}}$ is invertible, uniformly in $\s$.
As the $\e_{ij}\rightarrow 0$, these approximations improve and
$L_{g_{\sigma}}$
converges, as an operator on $L^2$, to the operator
$L_{g_{0}}$ on the disjoint union of
spheres indexed by $V$.

Writing the solution $g$  as a perturbation; $g=(1+\eta)^{\frac{4}{n-
2}}g_{\sigma}$ and using the conformal invariance of the conformal
Laplacian, we find that
$$
\aligned
L_{g}=&\,\Lap_{g} + n\\
=&\,(1+\eta)^{-\frac{4}{n-2}}\Lap_{g_{\sigma}} +
2(1+\eta)^{-\frac{n+2}{n-2}}\langle \nabla\eta,\nabla\,\rangle +n.
\endaligned  \tag 7.1
$$
Now, $\eta$ is small in $\Cal C^0$ (by which we mean that
$|\eta|_0$ can be
made arbitrarily small by taking all $\e_{ij}$ sufficiently small),
and estimates from [S2] can be used to show that $\nabla \eta$
is small in $L^2$ over any fixed compact set, but these
facts do not immediately imply the existence of an
analogous good basis for the  `small eigenspace' $K$ for $L_g$.
Nonetheless, we expect this basis to exist, namely that there
exist an $\e>0$,
such that if $\e_{ij}<\e$ for all $ij$ then there exists a
set of smooth functions, $\eta_{ij}$, $i \in V$, $j = 1, \cdots , n+1$,
which satisfy the estimates in Lemma 3.6 of [S2]. The span of these
functions is a subspace $K\subset L^{2}(\SL)$ such that
$$
\aligned
\|L_{g}\eta\|_{L^{2}(\SL)}\leq& c(\e)\|\eta\|_{L^{2}(\SL)} \text{ for }
\eta\in K\\
\|\eta\|_{L^{2}(\SL)}\leq& c\|L_{g}\eta\|_{L^{2}(\SL)} \text{ for }
\eta\in K^{\perp},  \endaligned
$$
where $c$ is a constant independent of $\e$, and $c(\e)$ tends to
$0$ as $\e \rightarrow 0$.
Note that from (7.1) this conjecture would be
immediate if $\|\eta\|_{H^{1}}$ were finite and small. Unfortunately,
this is never the case, since the $\eta$ in Schoen's
construction is not even in $L^2$.

The likeliest method to establish the existence of this space
$K$ with its explicit basis is to write the solution along
each end as a conformal perturbation of the Delaunay metric $g_{\e}$ on
the half-cylinder to which it converges at infinity.
There is good evidence that perturbation determining the solution
is uniformly small and even exponentially decaying,
again provided the neck sizes are sufficiently small.
With this information in hand, a not very complicated transference
procedure would produce the basis for $K$ with all necessary
estimates.

If these could be settled, it would then be possible to
resolve some basic questions about the
existence of an $L^2-$nullspace of $L_g$, for Schoen's solutions.
\proclaim{Conjecture 7.2} Suppose $g\in\ML$ is a solution similar
to that constructed in [S2], so that $g=(1+\eta)^{\frac{4}{n-2}}g_{\s}$
as above.
Assume that $\Lambda\subset\Sn$ does not lie in any
round hypersphere. Then there exists an $\e=\e(\Lambda)>0$, such that
if $\s$ satisfies $\e_{ij}<\e$ for each edge (i\.e\. all `neck sizes'
are less than $\e$) then $L_{g}$
has no $L^2-$nullspace, and hence $g$ is a smooth point of $\ML$.
If $\Lambda$ does lie in a hyperplane then any element $\phi$ in the
$L^2-$nullspace of $L_g$ is not integrable, i\.e\. it is not tangent to
a path in $\ML$.
\endproclaim

We provide a brief sketch of our plan for proving this. If $g$ is
a solution constructed from some $(\Gamma,\s)$, with all neck
sizes sufficiently small, and if $\phi \in L^2$, $L_g \phi = 0$,
then $\phi$ admits two different decompositions, one into
the `almost linear functions' $\eta_{ij}$ concentrated
on the almost spherical regions $\Omega_i$, and the other
(which is only local along each end $E_{\ell}$) into
eigenfunctions for the Laplacian along the cross-sectional
sphere, as in \S2. For this latter decomposition we need
the (putative) fact that the solution stays uniformly
close to the Delaunay solution to which it asymptotically
converges. The former expansion would show that $\phi$ is very
close to a linear function on each of the almost spherical
regions. A very important point here is that the linear
function determined by $\phi$ on each of the spherical
regions can never vanish. In terms of the eigenfunctions of
the second decomposition, this first decomposition implies
that along each end,
$\phi$ has most of its mass concentrated in the
eigenspaces corresponding to the zeroeth and first nonzero
eigenvalues of the cross-section.  Since $g$ is only
approximately Delaunay, the zeroeth Fourier coefficient
$\phi_0$ of $\phi$ is not an exact solution of $L_0\phi_0 = 0$,
however, the error terms are sufficiently controllable to
show that $\phi_0$ must be quite small (again, depending on
the size of $\e$). Hence, in fact, most of the mass of $\phi$
is contained in the coefficients $\phi_1, \cdots , \phi_{n}$,
corresponding the the first nonzero eigenvalue of the Laplacian
on $\SS^{n-1}$ (which has multiplicity $n$). We call these
`transverse' linear functions, since each of the spherical
regions where we consider them has a natural axis picked out.
Thus these arguments show that an $L^2$ solution
of $L\phi = 0$  restricts on each spherical region to be
almost linear. More specifically, on the innermost sphere
along each end (the one adjoining the central spherical region),
$\phi$ actually restricts to be approximately transverse linear.
But now, on the central sphere, $\phi$ is linear and near
each of the neighborhoods where the various ends are attached,
must restrict to a transverse linear function relative to
the axis determined by that end. This is clearly impossible
if the points do not lie in any subsphere.

For the second assertion, concerning the situation when $\Lambda
\subset \SS^{n-1}$, we only need to use that $\phi$ is
sufficiently close to a linear function on the central sphere.
Of necessity, this linear function is one which vanishes on the equator
determined by $\Lambda$. Hence, if $\phi$ were an integrable Jacobi
field, we would obtain a family of solutions which were not
symmetric about this equator, contradicting the reflectional
symmetry guaranteed by the Alexandrov reflection argument.

We hope to be able to settle these issues in the near future.
There are, of course, a number of more refined questions about
the singular structure of $\ML$, beginning with the problem,
already noted, of constructing solutions near which $\ML$ is
singular.

\medskip

\head Coordinates on $\ML$ \endhead

In this section we discuss two ways in which coordinates for $\ML$
may be given. We begin with the one arising from the linear analysis.

If $g \in \ML$ is a point for which $L_g$ has no $L^2-$nullspace,
then a neighborhood of $g$ in $\ML$ is parametrized by
a small neighborhood of the origin in the bounded nullspace $\Cal B$.
Thus, linear coordinates on $\Cal B$ yield local coordinates on $\ML$
near $g$. The problem then is to get precise information
about these linear coordinates on $\Cal B$. Unfortunately, this
seems difficult, in general. We now discuss how one would go about
this, and at least set the problem up in scattering-theoretic terms.

By definition, any element $\phi$
in the deficiency subspace $W$ can be expanded on each end as a
combination of the model solutions $\phi_1,\ \phi_2$ plus
an exponentially decreasing error. For clarity,
we label these model solutions on $E_j$ with a superscript $(j)$. Thus
$$
W \ni \phi \sim a_{j}\phi_1^{(j)} +  b_j\phi_2^{(j)}
\ \text{on\ } E_j. \tag 7.3
$$
So any $\phi \in W$ determines a map
$$
\gathered S: W \longrightarrow \RR^{2k} \\
 \phi \longmapsto (a_1, b_1, \dots , a_k, b_k) \endgathered \tag 7.4
$$
By the definition of $W$, $S$ is an isomorphism.
Now suppose $\phi$ is in the bounded nullspace $\Cal B$.
\proclaim{Lemma 7.5} Under the hypothesis (5.2), the linear map
$$
\left.C \equiv S \right|_{\Cal B}: \Cal B \longrightarrow \RR^{2k}
$$
is injective. Its image $S_{\Lambda}$
is a $k$-dimensional subspace which is Lagrangian with respect to the
natural symplectic form $\sum_{j=1}^k da_j \wedge db_j$ on $\RR^{2k}$.
\endproclaim
\demo{Proof} $C$ is obviously linear, and it is injective, because
otherwise there would exist an element $\phi \in \BL$ with
all its asymptotic coefficients $a_j, b_j$ vanishing. By the results
of the last section, such a $\phi$ would lie in $L^2$, which
contradicts the hypothesis (5.2).

That $S_{\Lambda}$ is Lagrangian is a consequence of Green's theorem.
There is an analogous result for manifolds with asymptotically cylindrical
ends proved in [Me]. Since we already know that $S_{\Lambda}$
is $k$-dimensional, it suffices to prove that it is isotropic.
For this, let $\psi, \ \chi \in \Cal B$ and set
$C(\psi) = (a_1, b_1, \dots, a_k, b_k),\ C(\chi) =
(\a_1, \b_1, \dots, \a_k, \b_k)$.
Using the coordinate $t_j$ along the end $E_j$, let $M_A$ denote the region
in $\SL$ where each $t_j \le A$. Then, since $L\psi = L\chi = 0$,
$$
0 = \int_{M_A} \left((L\psi) \chi - \psi(L\chi)\right)\,dV_g =
\int_{\del M_A} \left(\frac{\del \psi}{\del \nu}
\chi - \psi\frac{\del\chi}{\del \nu}\right)\,d\s_g,
$$
where $\nu$ is the unit vector-field normal to $\del M_A$. Substitute the
expansions for $\psi$ and $\chi$ at each end, and drop all terms
which decrease exponentially as $A \rightarrow \infty$ in this expression.
This means that we can replace the metric $g$ on $E_j$ by the Delaunay
metric $g_{\e_j}$ there, and also replace
$\nu$ by $u_{\e_j}^{-2/(n-2)}\del_{t_j}$ and $d\s_g$ by
$u_{\e_j}^{(2n-2)/n-2)}$. Taking
the limit as $A \rightarrow \infty$ we get
$$
0 = \sum_{j = 1}^k (\a_j b_j - a_j\b_j) \int_{\SS^{n-1}}
\left(\frac{\del \phi_1^{(j)}}{\del t_j}\phi_2^{(j)} - \phi_1^{(j)}
\frac{\del \phi_2^{(j)}}{\del t_j}\right)\,u_{\e_j}^2\,d\s. \tag 7.6
$$
Here $d\s$ is the standard volume measure on $\SS^{n-1}$. If we set
$A_j$ equal to the integrand for the $j^{\text{th}}$ end, then it
suffices to show that these constants $A_j$ are nonzero and independent
of $j$. However, $\phi_1^{(j)}$ and $\phi_2^{(j)}$ are
independent solutions of the ODE $L_0\phi = 0$, and the expression
in parentheses in the integrand in (7.6) is just the Wronskian for this
ODE. This Wronskian may be written explicitly as some nonzero multiple of
$$
\exp\left(-2\int \frac{u_{\e_j}'}{u_{\e_j}}\right) = u_{\e_j}^{-2}.
$$
This shows that the integrand, hence the integral itself, is a
nonvanishing constant, independent of $j$.
\enddemo

The relevance of $S_{\Lambda}$ to the description of
coordinates on $\ML$ near $g$ is explained as follows. Suppose
first, for simplicity, that there exists an element $\phi \in
\Cal B$ for which all $a_j,\,b_j = 0,\ j > 1$, but $a_1 = 1,\,b_1 = 0$.
If $g(s)$ is a curve in moduli space tangent to this $\phi$, then
the Delaunay parameters $\e_j$ on each of the ends would remain
fixed, at least infinitesimally, along this curve.
In fact, the metrics $g(s)$ (or the underlying curve of solutions
$u_s$ of the PDE) would be strongly asymptotic to the initial metric
$g$ on $E_j$, $j > 1$. On $E_1$ however, $g(s)$ would be strongly
asymptotic to a translation of $g$, or alternately, of its model
Delaunay metric there.  If, on the other hand, $a_1 = 0,\ b_1 = 1$ and
all other $a_j,\,b_j = 0$, then again each $g(s)$ would be strongly
asymptotic to $g$ on all ends except $E_1$; on that end the
Delaunay parameter for the model metric would
be changing, but the model Delaunay metric would not be translating.
Of course, neither of these types of elements of $\Cal B$ need exist.

If $\phi$ is a general element of $\Cal B$, and $g(s)$ is a curve in $\ML$
through $g$ tangent to $\phi$, then $\phi$ describes infinitesimal
changes in translation and Delaunay parameters along each end.
The most precise description of coordinates on $\ML$ would relate
the proportions of these various changes on each end to one another
for all directions $\phi \in T_g\ML$. This is equivalent to
describing the coefficients $a_j,b_j$ for each $\phi \in \Cal B$,
and this, in turn, is equivalent to describing the Lagrangian
subspace $S_{\Lambda}$. An explicit description of $S_{\Lambda}$
requires a better understanding of the Poho\v{z}aev balancing
condition (3.9).

For the sake of illustration, let us examine all this for $\Cal M_2$,
the moduli space when $k = 2$. Of course, we have already provided
a complete description of this space, but it provides a concrete
example of this description. We let the background metric be any
Delaunay metric $g_\e$ with $\e<\bar u$.
Now $\Cal B$ is $2$ dimensional, consisting
of the elements $\phi_1$ and $\phi_2$. Hence
$$
S_{\Lambda} = \{a_1 = a_2, b_1 = b_2\} \subset \RR^4,
$$
which is Lagrangian, as expected.
The two natural curves $g(s)$ emanating from any $g \in \Cal M_2$
are simply the families $\Phi_1$ and $\Phi_2$ defined in (2.14).

It may be possible to give a geometric description of parameters on $\ML$
by considering the Killing norms of the Poho\v{z}aev invariants.
Recall that these invariants are elements of $\frak o^{\star}(n+1,1)$,
the dual of the Lie algebra of the conformal group.  Any element
$X\in \frak o(n+1,1)$ may be uniquely decomposed as $X=X_{0}+{\text{\bf{w}}}$,
where $X_{0}\in \frak o(n+1)$ and ${\text{\bf{w}}}\in\RR^{n+1}$.
The Killing form
$$
B:\frak o(n+1,1)\times \frak o(n+1,1)\rightarrow \RR
$$
is the nondegenerate symmetric quadratic form given by
$$
B(X,\hat X)=\frac{1}{2} Tr(X_{0}\hat X_{0})+{\text{\bf{w}}}
\cdot{\hat{\text{\bf{w}}}}.
$$
Thus, $B(\cdot,\cdot)$ is positive definite on $\RR^{n+1}$ and
negative definite on $\frak o(n+1)$ (with respect to this decomposition).
Moreover, the adjoint representation preserves $B$ in the sense that
$B(Ad(F)(X),Ad(F)(\hat X))=B(X,\hat X)$,
for all $F\in O(n+1,1)$ and $X,\hat X\in \frak o(n+1,1)$.

Since $B$ is nondegenerate it provides an identification between
$\frak o(n+1,1)$ and its dual space.
Thus we may use the Killing form to identify, for the Poho\v{z}aev invariants
$\Cal P_{1},\dots,\Cal P_{k}$ of a metric $g\in\ML$, corresponding elements
$\Cal P^{\prime}_{1},\dots,\Cal P^{\prime}_{k}\in\frak o(n+1,1)$.
The squared killing norms of
the Poho\v{z}aev invariant $\Cal P_{i}$, is then defined to be
$$
\|\Cal P_{i}\|^{2} = B(\Cal P^{\prime}_{i},\Cal P^{\prime}_{i}).
$$
These $k$ numbers are conformal invariants.  By this we mean
the following.  If $F\in O(n+1,1)$ is a conformal diffeomorphism of $\Sn$ and
$g\in \ML$ then $F^{\star}(g)\in\Cal M_{F^{-1}(\Lambda)}$.  The Poho\v{z}aev
invariants of $g$ and $F^{\star}(g)$ do not coincide but transform via
the co-adjoint representation (see [S2] and [KKMS]).
Since $B$ is invariant under this representation this implies that
the squared killing norms of the Poho\v{z}aev invariants of $g$ and
$F^{\star}(g)$ are the same.

Schoen has suggested the following way to obtain parameters on $\ML$.
Let $\Lambda=\{p_1,\dots,p_k\}$ be a balanced
singular set. One can then try to produce, for some $\e>0$,
a 1-parameter family of solutions, $g_t\in\ML$, $t\in (0,\e)$,
such that the asymptotic neck sizes $\e_1,\dots,\e_k$ for $g_t$ are
all equal to $t$.  To realize the other $(k-1)$-parameters, choose any
$k$-tuple of numbers, $(a_1,\dots,a_k)$, close to 1 and normalized so that
$a_1=1$. Then there is a conformal diffeomorphism $F$ taking $\Lambda$
to $\Lambda^{\prime}=\{p^{\prime}_1,\dots,p^{\prime}_k\}$ so that
$\sum_{j=1}^{k}a_{j}p^{\prime}_j = 0$.
As before, there is now a 1-parameter family of solutions $g_{t}^{\prime}$,
on $\Sn\setminus\Lambda^{\prime}$ with asymptotic neck sizes given by
$(a_{1}t,\dots,a_{k}t)$. Thus $F^{\star}(g_{t}^{\prime})\in\ML$ is another
1-parameter family of solutions.  This exhibits the $k$-parameters as
$(t, a_2,\dots,a_k)$.  It should be possible to phrase this in terms of the
squared Killing norms of the Poho\v{z}aev invariants  described above.

\medskip

\head Dipole Solutions \endhead

We briefly describe a new construction of solutions [MPU]
for certain special configurations $\Lambda$.

Because the Delaunay solutions on the cylinder are nondegenerate,
in the sense that the linearized operators $L_{\e}$ have no
$L^2-$nullspace, it is relatively easy to graft them onto
known nondegenerate solutions. In this way, it is possible
to construct an explicit space of solutions when $\Lambda$
has $2k$ points which are grouped in pairs $\sigma_i$,
$i = 1, \dots k$, where the distance between points in each
pair is very much less than the distance between pairs.
We state the result when $k=2$ and the two pairs are assumed
to be located at opposite sides of $\Sn$.
\proclaim{Theorem 7.7 ([MPU])}  Given $\e>0$ there exists
an $\eta > 0$, such that for any $\Lambda = \{p_1,p_2,p_3,p_4\}$,
with $d(p_1,p_2) < \eta$,
$d(p_3,p_4) < \eta$, and $d(p_1,p_3) > 1$,
there exists a four parameter non-degenerate moduli
space $X_{\e}$ of solutions to the singular Yamabe
problem with neck-sizes at the points $p_i$ all
greater than $\e$. Moreover, the map from $X_{\e}$ to
$\RR^k$ sending a solution to the translation and
neck parameters of one point in each of the pairs,
say $p_1$ and $p_3$, is an isomorphism
onto an open set.
\endproclaim

These solutions can be described by making a conformal
change of the Delaunay solutions to $\RR^n \backslash
\{q_1,q_2\}$ and $\RR^n \backslash \{q_3,q_4\}$.
Simply adding the two solutions gives a result which is a
very close approximation to an exact solution which may be
constructed using the implicit function theorem.
It is necessary in this construction to make a small change in the
neck size and translation parameters of one point in each
pair to make the linear operators here surjective on
the relevant spaces.  Unfortunately, no uniqueness theorem
has been proved, even in this restricted setting.

It is of interest to note that
since the parameters on the other point of each pair may be
held fixed, this construction produces solutions with
asymptotically cylindrical ends.

\medskip
\head The Structure of the Boundary $\del \CML$ \endhead

There are a number of interesting questions about the structure
of the compactification $\CML$. Recall from \S3 that this
compactification is obtained by adding certain solutions
$g \in \Cal M_{\Lambda'}$ to $\ML$, for certain subsets
$\Lambda' \subset \Lambda$.  There are two basic problems.
The first is to determine which subsets $\Lambda'$ have
elements occurring in $\del \CML$. This is again intimately
connected to the balancing condition (3.9). There are few examples
where we can say anything explicit about this.
The case $k=2$ was already discussed in \S3.
The existence of dipole solutions [MPU] may give additional
information concerning this problem when $k$ is even.

The other fundamental problem is to determine how much of $\Cal M_{\Lambda'}$
is contained in $\del \CML$ when at least one point $g \in
\Cal M_{\Lambda'}$ lies in $\del \CML$. More specifically, we
propose the following:
\proclaim{Conjecture 7.8} Suppose that $g \in \Cal M_{\Lambda'} \cap
\del \CML$. Then the entire component of $\Cal M_{\Lambda'}$ containing
$g$ also lies in $\del \CML$.
\endproclaim
It is easy to
show, using the compactness theorem of [P2], that
the set of points $g'$, in the same component of $\Cal M_{\Lambda'}$
as $g$ and which also lie in $\del \CML$, is closed.
One method to prove that this set is open
would be to generalize the construction
of dipole solutions to a more general grafting procedure.
This would prove the conjecture.

\vfill\eject

\enddocument

\bye